\documentclass{article}

\usepackage{arxiv}
\usepackage{cite}
\usepackage{amsmath,amsfonts,bm}
\usepackage{amssymb}   
\usepackage{amsthm} 
\usepackage{mathtools}
\usepackage{algorithmic}
\usepackage{graphicx}
\usepackage{textcomp}
\usepackage{dsfont}
\usepackage{epstopdf}
\usepackage{orcidlink}
\usepackage{enumerate}
\usepackage{lscape}
\usepackage{multicol}
\usepackage{multirow}
\usepackage{booktabs}
\usepackage{framed}  
\usepackage{empheq}
\usepackage{graphics} % for pdf, bitmapped graphics files
\usepackage{epsfig} % for postscript graphics files
\usepackage{times} % assumes new font selection scheme installed
\usepackage{graphicx}
\usepackage{comment} 
\usepackage{apptools}
\usepackage{theoremref}
\usepackage{balance}
\graphicspath{{./},{./figures/}}

\newcommand{\tcb}{\textcolor{black}}
\newcommand{\trace}{{\operatorname{trace}}}
\newcommand{\factor}{\Phi}

\newtheorem{remark}{Remark}
\newtheorem{prop}{Proposition}

\newtheorem{cor}{Corollary}

\newtheorem{prob}{Problem}

\begin{document}

% \twocolumn[\begin{@twocolumnfalse}

\title{
The Holonomy of Optimal Mass Transport:\\
The Gaussian-Linear Case
}

\author{Mahmoud Abdelgalil\orcidlink{https://orcid.org/0000-0003-1932-5115}  and Tryphon T. Georgiou\orcidlink{https://orcid.org/0000-0003-0012-5447}
\thanks{Mahmoud Abdelgalil is with Electrical and Computer Engineering, University of California, San Diego, La Jolla, CA, USA, mabdelgalil@ucsd.edu.}
\thanks{Tryphon T. Georgiou is with Mechanical and Aerospace Engineering, University of California, Irvine, Irvine, CA, USA, tryphon@uci.edu.}}
\maketitle

\begin{abstract}
The theory of Monge-Kantorovich Optimal Mass Transport (OMT) has in recent years spurred a fast developing phase of research in stochastic control, control of ensemble systems, thermodynamics, data science, and several other fields in engineering and science.

We herein introduce a new type of transportation problems. The salient feature of these problems is that particles/agents in the ensemble are labeled and their relative position along their journey is of interest.
Of particular importance in our program are control laws that steer ensembles along cycles ensuring that individual particles return to their original position.
This feature is in contrast with the classical theory of optimal transport where the primary object of study is the path of probability densities, without any concern about particle labels.

In the theory that we present, we focus on the case Gaussian distributions and linear dynamics, and explore a hitherto unstudied sub-Riemannian structure of Monge-Kantorovich transport where the relative position of particles along their journey is modeled by the holonomy of the transportation schedule.
From this vantage point, we discuss several other problems of independent interest.
\end{abstract}
%
%%%%%%%%%%%%%%%%%%

\vspace*{.15in}
%%%%%%%%%%%%%%%%%%%
%

% \end{@twocolumnfalse}]

\section{Introduction}
The history of Optimal Mass Transport (OMT) can be traced to the work of Gaspar Monge in the latter part of the 18th century, who apparently was inspired by practical problems that he encountered in leveling dirt roads \cite{rachev1998monge,villani2003topics,villani2009optimal,chen2021stochastic}. It can be argued that Monge's problem constitutes one of the first optimal control problems, as he sought a transport schedule, for mass distributed in one location to be moved to another, that accrues minimal cost \cite{chen2021optimal}. Historically, Monge's cost was the actual Euclidean distance traversed by mass particles, although there are strong reasons to consider a (convex) cost, such as the square of the distance. Furthermore, Monge was interested in path planning besides the end to end correspondence. The history from there on was long, with the name of Kantorovich standing out --he invented duality theory to solve Monge's problem and received the {\em Nobel prize} for the transformative impact of his work.

The modern phase of OMT began in the 1990's with contributions by McCann, Gangbo, Benamou, Brenier, Otto, and many others \cite{mccann1997convexity,gangbo1996geometry,benamou2000computational,otto2001geometry}. 
\tcb{The control theoretic significance of the problem was soon recognized \cite{khesin2009nonholonomic,agrachev2009optimal,agrachev2009controllability,hindawi2011mass,rifford2014sub} and led to a flourishing research field that continues to this day.}
The form of the 
optimization problem that is pertinent herein, due to Benamou and Brenier \cite[Ch. 6]{villani2009optimal}, is:
\begin{subequations}\label{eq:minproblem12}
\begin{align}\label{eq:minproblem}
\min_{u(t,x)}\int_0^1 \int_{x\in\mathbb R^n} \mu_t(x) \|u(t,x)\|^2 dx dt,
\end{align}
subject to the continuity equation
\begin{align}\label{eq:minproblem2}
\partial_t \mu_t(x)=\nabla \cdot (\mu_t(x)u(t,x)),
\end{align}
\end{subequations}
and terminal conditions at $t\in\{0,1\}$; $n$ is the dimension of the space and $\nabla \cdot$ denotes the divergence.
Here, $u(t,\cdot)$ represents a vector field (our control action) that steers an initial probability density\footnote{Following a common slight abuse of notation one often uses the same symbol $\mu$ to denote the corresponding measure $\mu(dx)$, that in our case is always assumed absolutely continuous with respect to the Lebesgue.} $\mu_0(x)$ to a final one $\mu_1(x)$ over $t\in[0,1]$.

 The probability distributions that are being considered have finite second-moments, and the space of all such distributions (more generally, probability measures) is known as the Wasserstein space $\mathcal P_2(\mathbb R^n)$  \cite{villani2009optimal}.
The {\em square root} of the minimal value in \eqref{eq:minproblem} turns out to be a metric, the Wasserstein $\mathcal W_2(\mu_0,\mu_1)$ distance between $\mu_0$ and $\mu_1$. It coincides with
\begin{subequations}\label{eq:Monge}
\begin{align}\label{eq:W2}
\sqrt{\inf_{\varphi}\int_{\mathbb{R}^n}\lVert x-\varphi(x)\lVert^2 \,\mu_0(x)dx}
\end{align}
where the optimization takes place over maps
$
\varphi :\mathbb R^n\to \mathbb R^n
$
that ``push forward'' $\mu_0$ to $\mu_1$, in that
\begin{align}\label{eq:pushfwd}
\int_S\mu_1(x)dx= \int_{\varphi^{-1}(S)}\mu_0(x)dx,
\end{align}
\end{subequations}
for any set $S$ measurable. It is standard to denote this relation (push forward) by $\varphi \sharp \mu_0=\mu_1$. In fact, \eqref{eq:Monge} is the original formulation due to Monge.

The (unique) optimal transportation map\footnote{The optimal map exists under the assumption herein that $\mu_0,\mu_1$ are density functions, and thus, the corresponding measures are absolutely continuous with respect to the Lebesgue measure.} $\varphi^\star$ in \eqref{eq:W2} is known as the {\em Monge map} and it is of the form
\[
\varphi^\star(x)=\nabla \phi(x),
\]
where $\phi:\mathbb R^n\to \mathbb R$ is a convex function.
Moreover, the optimal path $\mu_t$ on $\mathcal P_2(\mathbb R^n)$ sought in \eqref{eq:minproblem12} that links the terminal specifications $\mu_0,\mu_1$, for $t\in[0,1]$, is
\begin{align}\label{eq:displacement_interpolation}
    \mu_t &= (\varphi_t){\sharp} \mu_0, \mbox{ where } \varphi_t= (1-t)\,\text{Id} + t \, \varphi^\star,
\end{align}
where $\text{Id}$ is the identity map. Thus, particles distributed according to $\mu_0(\cdot)$ are transported from their starting point $x$ to their final destination $\varphi^\star(x)$ along straight lines (i.e., along geodesics of the underlying space $\mathbb R^n$). In turn, the {\em McCann displacement interpolation} curves $\mu_t$ given in \eqref{eq:displacement_interpolation} turn out to be geodesics in their own right on $\mathcal P_2(\mathbb R^n)$ endowed with an (almost) Riemannian metric. This metric is due to Otto \cite{otto2001geometry,khesin2018geometric} and is
\begin{align}\label{eq:otto}
\langle \dot\mu_t,\dot\mu_t\rangle := \int_{x\in\mathbb R^n} \mu_t \|\nabla \phi_t\|^2 dx, 
\end{align}
where $\phi$ is the solution to the Poisson equation $\dot\mu_t=\nabla \cdot (\mu_t \nabla \phi_t)$. Note that $\dot\mu_t(\cdot)\equiv\partial_t\mu_t(\cdot)$ represents a tangent at $\mu_t$, and it directly corresponds to $\phi_t$ that can be identified with a  ``pressure distribution'' driving the flow \cite{otto2001geometry}, and to the optimal control  $u(t,x)=\nabla\phi_t(x)$ which is the force field that is generated by a time-varying control potential.

\begin{figure}[tp]
    \centering
    \includegraphics[width=0.4\linewidth]{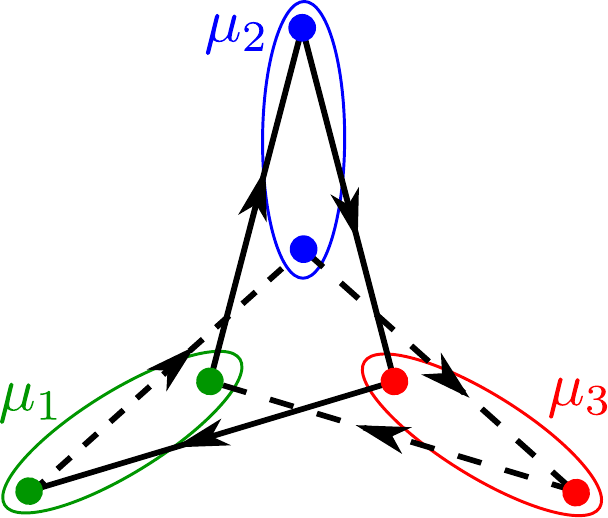} \hfill
    \caption{Trajectories of two tracer particles as the distribution consisting of two Dirac measures of equal weight traverses a triangle in $\mathcal{P}_2(\mathbb{R}^2)$ with sides corresponding to McCann geodesics.}
    \label{fig:mccann_triangle0}
\end{figure}

The starting point of our inquiry is the observation that when the distribution of an ensemble of particles traverses a closed orbit in $\mathcal P_2(\mathbb R^n)$ before it returns to its starting value, the particles do not necessarily return to their original positions. Indeed, traversing a closed path by an ensemble of particles, as a distribution, generically incurs a {  rearrangement of particle positions}. For instance, when considering a closed path of three McCann geodesics in $\mathcal P_2(\mathbb R^n)$ that connect pairs of $\bm{\mu}_1,\bm{\mu}_1,\bm{\mu}_3$ in succession, unless the Monge maps $\varphi_{i,j}^\star$ taking $\bm{\mu}_j$ to $\bm{\mu}_i$ commute, the composition
\[
\varphi:=\varphi_{3\mapsto 1}^\star\circ \varphi_{2\mapsto 3}^\star\circ \varphi_{1\mapsto 2}^\star
\]
will not be the identity. Thus, generically,  $\varphi(x)\neq x$, in spite of the fact that $\varphi\sharp \bm{\mu}_1=\bm{\mu}_1$, 
%$\varphi_{1,2}^\star\circ \varphi_{0,1}^\star\circ\varphi_{2,1}^\star\sharp \mu_1=\mu_1$ 
and similarly for a corresponding ordering in the case of $\bm{\mu}_2$ and $\bm{\mu}_3$.
The point is further exemplified in Fig.~\ref{fig:mccann_triangle0} by the rudimentary case of two ``tracer'' particles (equivalently, Dirac distributions) being transferred in succession between three specified configurations where they are in close proximity to one another, in a way that minimizes the overall Wasserstein distance being traversed.
It is observed that after one cycle their positions are exchanged. Clearly, in such cases where distributions consist of finitely many particles (singular distributions), closed orbits induce a permutation of the particles and it is the group of permutations that effects { rearrangement of particle positions}.
For more general probability measures on $\mathbb R^n$ or on a manifold, and under sufficient regularity conditions, the relevant groups are measure-preserving diffeomorphisms \cite{khesin2018geometric}.
{ Following fluid verbiage, we refer to the phenomenon outlined above, i.e., the rearrangement of particle positions without altering their distribution, as \emph{``mixing"}. The adopted nomenclature is not to be confused with the classical notion of mixing in ergodic theory where one studies the limiting behavior of repeated compositions of general measure preserving maps.
As will be shown in the sequel, the notion of mixing adopted herein turns out to have a natural interpretation as the \emph{holonomy} \cite{kobayashi1996foundations} of loops in the Wasserstein manifold $\mathcal P_2(\mathbb R^n)$.}
The present paper is devoted to the study of said holonomy in the case of Gaussian probability distributions. That is, we study optimal transport in the fibre bundle of {  linear} transportation plans, where a holonomy group acts naturally, and focus on developing a suitable differential geometric framework to { study such rearrangement of particle positions}. 
Specifically, in Section \ref{sec:II} we provide a brief expository on optimal transport of Gaussian distributions highlighting the holonomy that arises when traversing a closed path in the Wasserstein space. In Section \ref{sec:III} we develop the differential geometric setting that is suitable for studying parallel transport and holonomy along Wasserstein curves. In Section \ref{sec:IV} we detail { a hitherto unstudied sub-Riemannian structure induced by optimal transport wherein the problem of realizing a specified parallel transport map along optimal curves in the Wasserstein space can be naturally formulated and studied.} Section \ref{sec:Applications} overviews examples and applications that highlight the relevance of the theory we presented. We close with an epilogue, in Section \ref{sec:Epilogue}, where we envision enabling new aspects of the framework and concepts discussed in the paper.

\section{Transport in the space of Gaussian distributions}\label{sec:II}

In the present work we focus on the transport of distributions by way of potential forces. That is, the basic model that we consider at the level of transport of individual particles is a linear system
\begin{align}\label{eq:system}
    \dot X_t &= -\nabla \phi_t(X_t) =: u_t\\
    & = A_t X_t + r_t,\nonumber
\end{align}
with $X_t\in \mathbb R^n$ distributed according to the Gaussian distribution $\mu_t$ (denoted by $X_t\sim\mu_t$) having covariance $\Sigma_t$ and mean $m_t$. The dynamics generated by a quadratic time-varying control potential $\phi_t(x)=-x^\intercal A_t x + x^\intercal r_t$.
It is noted that in \eqref{eq:system} there is no added stochastic excitation, and the particles obey deterministic dynamics, in this Lagrangian representation.

Because the mean is directly regulated by the reference $r_t$,
as
$
\dot m_t=A_t m_t +r_t,
$
we only focus on the flow of covariances. That is, we assume throughout that $r_t,m_t$ vanish and that 
\begin{align}
\phi_t(x)=-x^\intercal A_t x. \label{eq:potential}   
\end{align}
Thus, Gaussian distributions have zero mean and covariances $\Sigma_t=\mathbb E\{X_tX_t^\intercal\}$ that obey the Lyapunov equation
\begin{align}\label{eq:Lyapunov}
    \dot \Sigma_t&=A_t\Sigma_t + \Sigma_t A_t^\intercal.
\end{align}
Moreover, $A_t$ will be taken to be symmetric, i.e., $A_t=A_t^\intercal$. Equation \eqref{eq:Lyapunov} represents an {\em Eulerian specification} since it does not encode the actual particle location. The flow of particles is effected by a state transition matrix $\factor_t$, where 
\begin{align}\label{eq:statetransition}
\dot\factor_t=A_t\factor_t
\end{align}
and $\Phi_0=I$, though occasionaly we use a double index for
$\dot\factor_{t,s}=A_t\factor_{t,s}$ for $\factor_{s,s}=I$. Evidently, $\Sigma_t=\factor_{t,s} \Sigma_s\factor_{t,s}^\intercal$.

Thus, henceforth, we will consider non-degenerate Gaussian distributions {\em centered} at the origin. We will denote by $\mathcal{N}(\mathbb{R}^n)$ the space of all such distributions in $\mathbb{R}^n$, i.e., of the form
\begin{align}\label{eq:gaussian_distribution}
    \mu(x) &= \frac{\exp\left(-\frac{1}{2} x^\intercal \Sigma^{-1} x\right)}{\sqrt{(2\pi)^n\text{det}(\Sigma)}},
\end{align}
with $\Sigma \in\text{Sym}^+(n)$,  the symmetric and positive definite matrices
\begin{align}\nonumber
    \text{Sym}^+(n):=\{\Sigma\in\mathbb{R}^{n\times n}~|~\Sigma=\Sigma^\intercal,\,\Sigma\succ 0\}.
\end{align}

\subsection{Optimal transport of Gaussian distributions}

The problem to relate optimally two Gaussian distributions has received attention by several authors \cite{knott1984optimal,olkin1982distance,givens1984class,dowson1982frechet}, and the solution to Monge's OMT in the setting of Gaussian distributions has arisen from several different angles. For instance, Fr\'echet's 
distance between probability distributions
\[
\min_{\Pi}\mathbb E\{\|X_0-X_1\|^2\}
\]
where optimization is over $\Pi(dx_0,dx_1)$, the joint probability law of two random variables with given marginal distributions, $X_0\sim \mu_0$ and $X_1\sim \mu_1$, can be seen as a special case of the regularization to Monge's problem \eqref{eq:Monge} introduced by Kantorovich \cite{kantorovich2006translocation}. When the cost is quadratic, whether in the form
\eqref{eq:Monge} or \eqref{eq:minproblem12}, the solutions coincide. This has been explained in the context of the more general question of interpolation in $\mathcal P_2(\mathbb  R^n)$, see \cite{mccann1997convexity,villani2003topics,villani2009optimal}.

Specifically, for two (zero mean) Gaussian distributions $\mu_0,\mu_1\in\mathcal{N}(\mathbb{R}^n)$ with covariances $\Sigma_0,\Sigma_1\in\text{Sym}^+(n)$, respectively,
the optimal transportation map $\varphi^\star$ in Monge's problem \eqref{eq:Monge} is linear,
\begin{align}
    \varphi^\star(x)&= \factor^\star x
\end{align}
with $\factor^\star$ the symmetric positive definite matrix\footnote{The
\label{footnote3}
two expressions \eqref{eq:optimalMonge} and \eqref{eq:optimalMonge2} are equivalent and characterize the maximal $\Phi=\Phi^\intercal$, in the positive definite sense, that satisfies 
\[\begin{bmatrix}
    M_0 & \Phi\\\Phi & M_1
\end{bmatrix}\geq 0,\]
for $M_0=\Sigma_0^{-1}$ and $M_1=\Sigma_1$.
This symmetric $\Phi$ is also referred to as the  {\em geometric mean} of the nonnegative matrices $M_0$ and $M_1$, see \cite{ando2004geometric} for references and generalizations to multiple matrices.}
\begin{subequations}
\begin{align}\label{eq:optimalMonge}
    \factor^\star&= 
    \Sigma_{0}^{-\frac{1}{2}}
    \big(\Sigma_{0}^{\frac{1}{2}}\Sigma_1\Sigma_{0}^{\frac{1}{2}}\big)^{\frac12}
    \Sigma_{0}^{-\frac{1}{2}}\\\label{eq:optimalMonge2}
    &=\Sigma_{1}^{\frac{1}{2}}
    \big(\Sigma_{1}^{\frac{1}{2}}\Sigma_0\Sigma_{1}^{\frac{1}{2}}\big)^{-\frac12}
    \Sigma_{1}^{\frac{1}{2}},
\end{align}
giving
\[
\mathcal W_2(\mu_0,\mu_1)=\trace(\Sigma_0+\Sigma_1-2(\Sigma_0^\frac12\Sigma_1\Sigma_0^\frac12)^\frac12).\]
\end{subequations}

In light of \eqref{eq:displacement_interpolation}, $\varphi_t(x)= ((1-t)I+\factor^\star)x$,  and McCann's displacement interpolation  $(\varphi_t){\sharp} \mu_0$ remains in  $\mathcal{N}(\mathbb{R}^n)$ for $t\in[0,1]$. Specifically, this is
\begin{align*}
     \mu_t(x) &= \frac{\exp\left(-\frac{1}{2} x^\intercal \Sigma_t^{-1} x\right)}{\sqrt{(2\pi)^n\text{det}(\Sigma_t)}},
\end{align*}
with covariance
\begin{align}\label{eq:Sigma_geodesic}
    \Sigma_t&= ((1-t) I + t \factor^\star){\sharp} \Sigma_0,
\end{align}
where $\Phi{\sharp}\Sigma$ denotes matrix congruence,
\begin{align*}
    \Phi{\sharp}\Sigma = \Phi \Sigma \Phi^\intercal,
\end{align*}
echoing the correspondence in $\varphi \sharp \mu$.
Finally, from \eqref{eq:Sigma_geodesic}, $\Sigma_t$ obeys \eqref{eq:Lyapunov} for $\in[0,1]$ and a suitable $A_t$, namely 
\begin{align}\label{eq:A}
A_t=(\factor^\star-I)(I+t(\factor^\star-I))^{-1}.
\end{align}
 We see that the optimal control in \eqref{eq:minproblem} is $u(t,x)=A_t x$ while \eqref{eq:minproblem2} reduces to \eqref{eq:Lyapunov}.

\begin{remark}\label{rem:remark1}
    It can be seen from \eqref{eq:A} and \eqref{eq:optimalMonge} that, provided $\factor^\star-I$ is not singular, $A_t$ is invertible in the interval $t\in[0,1]$, that is of interest. In fact, $A_t$ remains invertible outside $[0,1]$, as explained next, and the McCann geodesics (i.e., $\Sigma_t$) can be extrapolated accordingly.
    Specifically, if $\Sigma_0<\Sigma_1$, then $A_t$ remains invertible for all $t\in[0,\infty)$, otherwise, it remains invertible as long as $t<(1-\lambda_{\rm max}^{-\frac12})^{-1}$, with $\lambda_{\rm max}$ the largest eigenvalue of $\Sigma_0\Sigma_1^{-1}$.
Similarly, when $\Sigma_0>\Sigma_1$, $A_t$ is invertible for all $t\in(-\infty,0]$,  otherwise it remains invertible as long as $t>(\lambda_{\rm min}^{-1/2}-1)^{-1}$, where $\lambda_{\rm min}$ is the smallest eigenvalue of $\Sigma_0\Sigma_1^{-1}$. Evidently, there is no flow along directions where $\factor^\star x=x$, and hence, this trivial case does not occur when $\factor^\star-I$ is invertible. $\Box$
\end{remark}

\subsection{The Holonomy of Gaussian triangles}

We recapitulate the starting point of our inquiry, from the introduction, by exploring triangles in $\mathcal{N}(\mathbb{R}^n)$ with edges constructed via McCann interpolation between vertices.

Identifying distributions in $\mathcal{N}(\mathbb{R}^n)$ with their corresponding covariances, we begin by taking vertices $\mathbf \Sigma_1,\mathbf\Sigma_2,\mathbf\Sigma_3\subset\text{Sym}^+(n)$ and construct edges between pairs of $\mathbf\Sigma$'s using McCann geodesics. For brevity, we refer to a closed cycle  $\Sigma_t$ traversing the edges over $t\in[0,1]$ as a {\em Gaussian triangle}.
The time-indexing in traversing edges is immaterial, so we choose equal ``time steps'' for traversing each edge and set
% This, we denote $\triangle:[0,1]\rightarrow \text{Sym}_+(n)$ defined by
\newcommand{\ab}{{1\mapsto 2}}
\newcommand{\bc}{{2\mapsto 3}}
\newcommand{\ca}{{3\mapsto 1}}
\begin{align}
    \Sigma_t&:=\begin{cases}
        (I + 3t (\factor_\ab^\star-I)){\sharp}\mathbf\Sigma_1, & t\in\left[0,\frac{1}{3}\right),\\
        (I + 3(t-\frac{1}{3}) (\factor_\bc^\star-I)){\sharp}\mathbf\Sigma_2, & t\in\left[\frac{1}{3},\frac{2}{3}\right),\\
        (I + 3(t-\frac{2}{3}) (\factor_\ca^\star-I)){\sharp}\mathbf\Sigma_3, & t\in\left[\frac{2}{3},1\right],\\
    \end{cases}
\end{align}
with $\factor_\ab^\star,\factor_\bc^\star,\factor_\ca^\star$ the Monge maps between vertices.

\begin{figure}[t]
    \centering
    \includegraphics[width=0.4\linewidth]{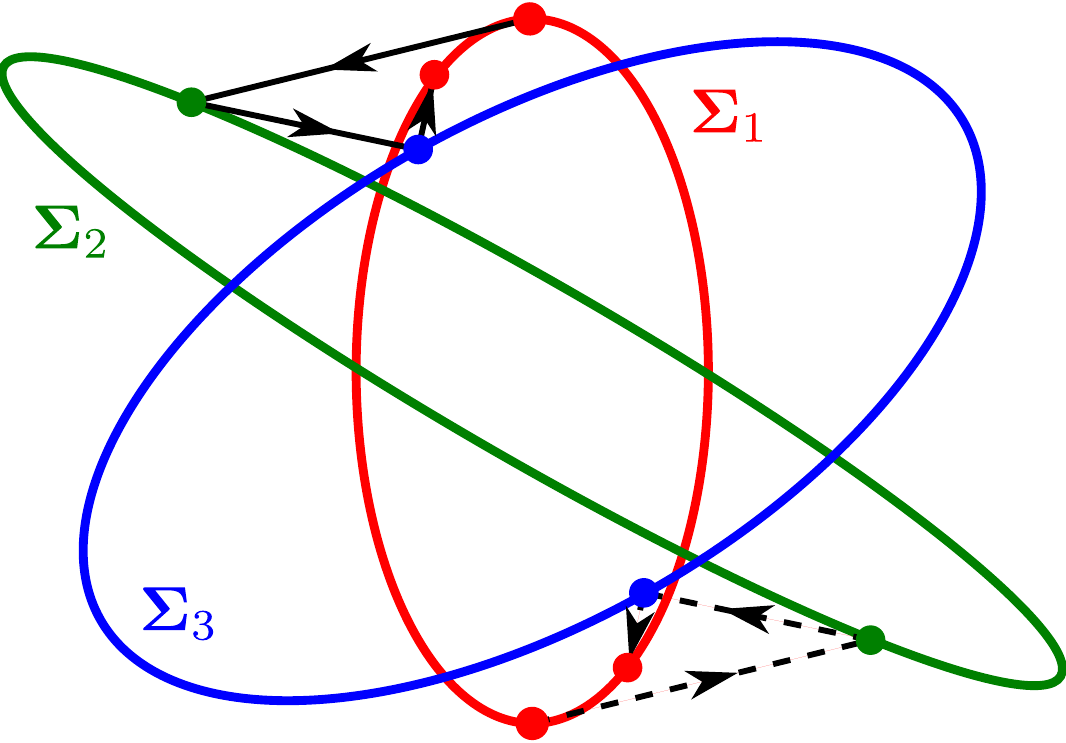} 
    \caption{Trajectories of two tracer particle traversing a Gaussian triangle with vertices $\{\mathbf\Sigma_1,\mathbf\Sigma_2,\mathbf\Sigma_3\}$. The arrows indicate the direction of motion of the tracer particles. Each colored ellipse represents a level set of the corresponding Gaussian probability density.}
    \label{fig:mccann_triangle}
\end{figure}

After traversing the cycle, $\Sigma_1=\Sigma_0 =\mathbf\Sigma_1$. The composition
\begin{align*}
\Theta&=\factor_\ca^\star\factor_\bc^\star\factor_\ab^\star,
\end{align*}
is a state transition matrix that satifies
\begin{align}\label{eq:SigmaUnitary}
\Theta\sharp \mathbf\Sigma_1=\mathbf\Sigma_1,
\end{align}
as can be readily verified using \eqref{eq:optimalMonge}. But $\Theta\neq I$, in general. Thus, $X_1\neq \Theta X_0$ for a tracer particle starting at $X_0=x$.
% Suppose further that we track the position of a \emph{tracer} particle as the distribution of particles traverses the McCann triangle once, starting at $\Sigma_1$. If the particle's initial position is $x_{0}$, then its trajectory is
% \begin{align*}
%     x_t&:=\begin{cases}
%         ((1-3t) I + 3t \factor_{1\mapsto 2}^\star) x, & t\in\left[0,\frac{1}{3}\right),\\
%         ((2-3t) I + (3t-1) \factor_{2\mapsto 3}^\star)\factor_{1\mapsto 2}^\star x, & t\in\left[\frac{1}{3},\frac{2}{3}\right),\\
%         ((3-3t) I + (3t-2) \factor_{3\mapsto 1}^\star)\factor_{2\mapsto 3}^\star\factor_{12}^\star x, & t\in\left[\frac{2}{3},1\right].\\
%     \end{cases}
% \end{align*}
% In particular, the final position of the particle is given by
% \begin{align}
%     x_1 &= \Theta x_0, & \Theta&=\factor_{3\mapsto 1}^\star\factor_{2\mapsto 3}^\star\factor_{1\mapsto 2}^\star
% \end{align}
% Clearly, $x_1\neq x_0$ in general despite the fact that
% \begin{align}
%     \Theta{\sharp} \Sigma_1&=\Sigma_1,
% \end{align}
% as can be readily verified. 
The two-dimensional sketch of Figure \ref{fig:mccann_triangle} highlights the point.

Obviously, mixing is not exclusive to triangles. A succession of state transition maps that allows traversing a polygon, or any closed curve
$\Sigma_t\in\text{Sym}^+(n)$
for that matter, i.e.,
\[
\Sigma_1=\factor_{01}\Sigma_0\factor_{01}^\intercal=\Sigma_0,
\]
leaves $\Sigma_0$ invariant under congruence. Such state transition matrices are special and will typically be denoted using the symbol $\Theta$. Hence, here, we write $\Theta=\factor_{01}$ and observe that
\[
\Sigma_0^{-\frac12}\Theta\Sigma_0^{\frac12}
\]
is an orthogonal matrix. In fact, $\Theta$ has positive determinant, and hence equal to $1$, since it arises by integrating \eqref{eq:statetransition} from the identity.
% is any geodesic polygon, i.e. a closed curve that connects a collection of covariances $\{\Sigma_i\}_{i=1}^N$ in the order 
% $$\Sigma_1\rightarrow\Sigma_2\rightarrow \cdots\rightarrow\Sigma_N\rightarrow \Sigma_1,$$ 
% such that the segment between any two successive points is a McCann geodesic, then the displacement of a tracer particle from its initial position after traversing the polygon is given by $x_1= \Theta x_0$ where $\Theta$ is the matrix
% \begin{align}\label{eq:example_measure_preserving_map}
%     \Theta &= \factor_{N\mapsto 1}\prod_{i=1}^{N-1}\factor_{i\mapsto i+1}.
% \end{align}
%A matrix $\Theta$ that satisfies $\Theta{\sharp} \Sigma=\Sigma$ for some $\Sigma \in \text{Sym}_+(n)$ is associated with a (linear) 
It corresponds to a
\emph{measure-preserving} map $\varphi(x)=\Theta x$ that preserves the Gaussian distribution with covariance $\Sigma_0$ under the pushforward operation. 
The set of all such matrices constitutes a (closed and connected) Lie subgroup of the general linear group $\text{GL}(n)$, this is the group
\begin{align}
    \text{SO}(n,\Sigma):=\{\Theta\in \text{GL}(n)~|~\Theta{\sharp}\Sigma = \Sigma,\; \det(\Theta)=1\}. 
\end{align}

%By definition, any matrix of the form \eqref{eq:example_measure_preserving_map} is an element of the group $\text{O}(n,\Sigma_1)$. 

A natural question that arises in the present context is whether it is possible to realize any $\Theta\in \text{SO}(n,\Sigma)$ in this way, solving \eqref{eq:statetransition} via a choice of a control protocol $A_t$, for $t\in[0,1]$, possibly piecewise constant.
This question is equivalent to the classical control-theoretic question of controllability for a certain driftless control-affine system. Finding an ``optimal" closed curve that achieves a given $\Theta\in \text{SO}(n,\Sigma)$ reduces to finding a \emph{sub-Riemannian geodesic} \cite{agrachev2019comprehensive} or, equivalently, the solution to an \emph{isoholonomic problem} in a certain principal bundle \cite{montgomery1990isoholonomic,montgomery2002tour}. Therefore, the natural framework wherein one can formulate and study mixing in OMT is the \emph{sub-Riemannian geometry}  \cite{agrachev2019comprehensive, montgomery2002tour} of \emph{principal bundles} \cite{sontz2015principal,kobayashi1996foundations}. The next section gives an overview of the requisite mathematical machinery, specialized to the setting of covariance matrices and linear transformations. To enhance readability, we opt to present our results in coordinates rather than in a coordinate-free fashion as is customary in differential geometry. Our choice is justified by the fact that the manifolds we work with have global charts. 

\section{Differential Geometry of Linear Transport}\label{sec:III}

 In \cite{otto2001geometry}, Otto showed that the Wasserstein metric $W_2$ is (formally) induced by a Riemannian metric, subsequently referred to as the Wasserstein-Otto metric \cite{khesin2018geometric} or the Wasserstein metric \cite{takatsu2011wasserstein}, which can be obtained via a \textit{Riemannian submersion} of a flat Riemannian metric on $\text{Diff}(\mathbb{R}^n)$, the diffeomorphism group of $\mathbb{R}^n$. His insight led to a vibrant body of literature on the differential geometry of OMT, Riemannian geometry of Wasserstein spaces, and geometric hydrodynamics. Otto's geometric framework is a natural avenue for our purposes. For that reason, we now give an overview of the pertinent components of Otto's framework, specialized to Gaussian distributions and linear diffeomorphisms.

% \subsection{Gaussian distributions}
As noted earlier, we may identify $\mathcal{N}(\mathbb{R}^n)$ with $\text{Sym}^+(n)$, since 
  $\mathcal{N}(\mathbb{R}^n)$ is globally parameterized by $\text{Sym}^+(n)$ from \eqref{eq:gaussian_distribution}. 
But $\text{Sym}^+(n)$ is an open subset of the vector space of symmetric matrices
\begin{align}
    \text{Sym}(n):=\{S\in\mathbb{R}^{n\times n}~|~S = S^\intercal\}.
\end{align}
Therefore, $\text{Sym}^+(n)$ is a smooth manifold and, at any $\Sigma\in\text{Sym}^+(n)$, its tangent space
$T_\Sigma \text{Sym}^+(n)$
 is canonically isomorphic to $\text{Sym}(n)$ \cite{lee2012smooth}. 
 
 The group of linear diffeomorphisms of $\mathbb R^n$ is the finite dimensional subgroup of $\text{Diff}(\mathbb{R}^n)$, namely,
\begin{align*}
    \text{LDiff}(\mathbb{R}^n):=\{\varphi\in\text{Diff}(\mathbb{R}^n)~|~\varphi(x) = \factor x, \,\factor\in\text{GL}(n)\}.
\end{align*}
%where $\text{GL}(n)$ is the general linear group. 
It has two connected components, with the identity component being
\begin{align*}
    \text{LDiff}^+(\mathbb{R}^n):=\{\varphi\in\text{Diff}(\mathbb{R}^n)~|~\varphi(x) = \factor x, \,\factor\in\text{GL}^+(n)\},
\end{align*}
where as usual,
$
    \text{GL}^+(n):=\{\factor\in\text{GL}(n)~|~\det(\factor) > 0\},
$
is the identity component of $\text{GL}(n)$ (i.e., matrices with positive determinants).
An element $\varphi\in\text{LDiff}^+(\mathbb{R}^n)$ acts on $\mathcal{N}(\mathbb{R}^n)$ from the left via the pushforward on $\mathcal{N}(\mathbb{R}^n)$, namely,
\begin{align}
    L_\varphi: \mu \mapsto \varphi{\sharp}\mu,
\end{align}
It is clear that $\text{LDiff}^+(\mathbb{R}^n)$ is globally parameterized by $\text{GL}^+(n)$, which is a Lie group.
\begin{subequations}\label{eq:left_action_pushforward}
Moreover, direct computation of the pushforward $\varphi{\sharp}\mu$ yields
\begin{align}
    \varphi{\sharp}\mu(x) = \frac{\exp\left(-\frac{1}{2} x^\intercal ( \factor\sharp \Sigma)^{-1} x\right)}{\sqrt{(2\pi)^n\text{det}(\factor\sharp\Sigma)}} \in \mathcal{N}(\mathbb{R}^n),
\end{align}
when $\varphi\in\text{LDiff}^+(\mathbb{R}^n)$ and $\mu \in\mathcal{N}(\mathbb{R}^n)$ are
\begin{align}
    \varphi(x)&=\factor x, & \mu(x) &= \frac{\exp\left(-\frac{1}{2} x^\intercal \Sigma^{-1} x\right)}{\sqrt{(2\pi)^n\text{det}(\Sigma)}}.
\end{align}
\end{subequations}

In conclusion, from \eqref{eq:left_action_pushforward}, we observe that the left action of $\varphi\in\text{LDiff}^+(\mathbb{R}^n)$ on $\mathcal{N}(\mathbb{R}^n)$ via the pushforward is equivalent to the left action of $\factor\in\text{GL}^+(n)$ on $\text{Sym}^+(n)$ via congruence:
\begin{align}\label{eq:left_action_congruence}
    L_\factor:
    %\text{Sym}^+(n)\ni 
    \Sigma \mapsto \factor\sharp\Sigma.
    %\in \text{Sym}^+(n).
\end{align}
With that equivalence in mind, we henceforth work with the parameterizations $\text{GL}^+(n)$ and $\text{Sym}^+(n)$ instead of $\text{LDiff}^+(\mathbb{R}^n)$ and $\mathcal{N}(\mathbb{R}^n)$, respectively. 

\subsection{The OMT principal bundle}

\vspace{0.2cm}
We begin by fixing an arbitrary reference covariance $\Sigma_{\text{ref}}\in \text{Sym}^+(n)$. We then define the map
\begin{align}
    \pi\;:\;\text{GL}^+(n)\to \text{Sym}^+(n)\;:\;
    \factor\mapsto L_\factor(\Sigma_{\text{ref}}).
\end{align}
One may take $\Sigma_{\text{ref}}$ to be the identity, and then $\pi\;:\;\factor\mapsto \factor\factor^\intercal$. The map $\pi$ will serve to define the principle bundle structure with base space the set of covariances, $\Sigma=\factor\factor^\intercal$, or $\Sigma=\factor\Sigma_{\text{ref}}\factor^\intercal$ with respect to any reference $\Sigma_{\text{ref}}$.

The map $\pi$ is in fact a restriction of the map introduced by Otto in his seminal paper \cite{otto2001geometry} to our finite dimensional setting. In particular, $\pi$ is a smooth submersion, and so $\pi:\text{GL}^+(n)\rightarrow\text{Sym}^+(n)$ is a smooth fiber bundle. For any $\Sigma\in\text{Sym}^+(n)$, the fiber at $\Sigma$ is the submanifold:
\begin{align}
    \pi^{-1}(\Sigma) = \{\factor\in\text{GL}^+(n)~|~\pi(\factor) = L_{\factor}(\Sigma_{\text{ref}}) = \Sigma\}.
\end{align}
In particular, the fiber $\pi^{-1}(\Sigma_{\text{ref}})$ at the reference is the \textit{isotropy group} of $\Sigma_{\text{ref}}$ and, evidently, coincides with 
\begin{align}
    \text{SO}(n,\Sigma_{\text{ref}})=\{\Theta\in\text{GL}^+(n)~|~ L_\Theta(\Sigma_{\text{ref}}) = \Sigma_{\text{ref}}\},
\end{align}
repeated here for clarity,
%$\text{SO}(n,\Sigma_{\text{ref}})$, 
under the left action of $\text{GL}^+(n)$.
As a (closed) subgroup of $\text{GL}^+(n)$, $\text{SO}(n,\Sigma_{\text{ref}})$ acts naturally from the right on $\text{GL}^+(n)$ via matrix multiplication:
\begin{align}
    R_\Theta\;:\;\text{GL}^+(n)\to \text{GL}^+(n)\;:\; \factor \mapsto \factor\Theta.
\end{align}
This action is \textit{free} since $\text{SO}(n,\Sigma_{\text{ref}})$ is a subgroup of $\text{GL}^+(n)$. In addition, if $\Sigma\in \text{Sym}^+(n)$, $\factor\in\pi^{-1}(\Sigma)$, and $\Theta\in \text{SO}(n,\Sigma_{\text{ref}})$, then 
\begin{align}
    \pi\circ R_\Theta(\factor) = \factor\Theta \Sigma_{\text{ref}} \Theta^\intercal \factor^\intercal = \factor\Sigma_{\text{ref}} \factor^\intercal =  \Sigma,
\end{align}
which implies that $R_\Theta(\factor)\in\pi^{-1}(\Sigma)$. That is, the right action of $\text{SO}(n,\Sigma_{\text{ref}})$ on $\text{GL}^+(n)$ preserves the fibers of the bundle $\pi:\text{GL}^+(n)\rightarrow\text{Sym}^+(n)$.

Finally, it can be shown that the action of $\text{SO}(n,\Sigma_{\text{ref}})$ on each fiber is \emph{transitive} \cite{modin2016geometry}, i.e. if $\factor$ and $\tilde{\factor}$ belong to the same fiber $\pi^{-1}(\Sigma)$, then there exists $\Theta\in \text{SO}(n,\Sigma_{\text{ref}})$ such that $\tilde{\factor}=\factor\Theta$. Since the action of $\text{SO}(n,\Sigma_{\text{ref}})$ on the fibers is both free and transitive, the following proposition is immediate \cite{montgomery2002tour,kobayashi1996foundations}.
\vspace{0.1cm}
\begin{prop}
    The fiber bundle $\pi:\text{GL}^+(n)\rightarrow\text{Sym}^+(n)$ is a principal $\text{SO}(n,\Sigma_{\text{ref}})$-bundle over $\text{Sym}^+(n)$.
\end{prop}
\vspace{0.1cm}
Henceforth, with a slight abuse of notation, we use $\pi$ to refer both the map $\pi$ and to the fiber bundle $\pi:\text{GL}^+(n)\rightarrow\text{Sym}^+(n)$. The meaning will be clear from the context.
\begin{remark}
The statement of the proposition, for the special case where $\Sigma_{\text{ref}}=I$, amounts to the statement that matrices $\Phi$ with positive determinant constitute a principal bundle of factors over the base space of covariances $\Sigma=\Phi\Phi^\intercal$, and that the holonomy group of the fibers is a subgroup of the special orthogonal group of matrices $\Theta$, since $\factor\to \factor\Theta$ leaves $\pi(\factor\Theta)$ the same and equal to $\Sigma$. \hfill $\Box$
\end{remark}

\subsection{The OMT Ehresmann connection}

We now consider the differential $\text{d}\pi$ of $\pi$ that maps the tangent space of 
$\text{GL}^+(n)$ to that of $\text{Sym}^+(n)$. Our interest is in identifying permissible directions for elements $\Phi$ on the fibers $\pi^{-1}(\Sigma)$ above covariances $\Sigma$ on the base space. Thus, it is essential to relate elements on nearby fibers and how one is allowed to transition between those, and this is done with the concept of an Ehresmann connection.

The kernel of $\text{d}\pi$ defines the so-called \textit{vertical sub-bundle}
\begin{align*}
    \text{Ver}_\factor:= \text{Ker}(\text{d}\pi_\factor) = \{\dot{\factor}\in T_\factor\text{GL}^+(n)~|~\text{d}\pi_\factor(\dot{\factor}) = 0\}.
\end{align*}
A \textit{principal Ehresmann connection} \cite{kobayashi1996foundations,sontz2015principal}, or \emph{principal connection} for short, allows transporting between different fibers and is based on a choice of a {\em horizontal sub-bundle}.
Specifically, this is a sub-bundle $\text{Hor}$ that associates to every $\factor\in\text{GL}^+(n)$ a subspace $\text{Hor}_\factor\subset T_\factor\text{GL}^+(n)$ so that, for all $\factor\in\text{GL}^+(n)$ and all $\Theta\in\text{SO}(n,\Sigma_{\text{ref}})$, the following hold:
\begin{enumerate}
    \item $T_\factor\text{GL}^+(n) = \text{Ver}_\factor \oplus \text{Hor}_\factor$,
    \item $\text{Hor}_{\factor\Theta}=(R_\Theta)_*\text{Hor}_{\factor}$.
\end{enumerate}
Intuitively, the last condition can be interpreted as a property of the horizontal spaces at different points on the fiber being parallel, in a suitable sense.

In words, a principal Ehresmann connection is a \emph{canonical} choice of a sub-bundle complementary to the vertical bundle $\text{Ver}:\text{GL}^+(n)\ni \factor\mapsto \text{Ver}_{\factor}\subset T_{\factor}\text{GL}^+(n)$ and invariant under the action of the isotropy group $\text{SO}(n,\Sigma_{\text{ref}})$. There is considerable freedom in the construction of the sub-bundle defining the Ehresmann connection, i.e., the \textit{horizontal sub-bundle}. However, here, we are primarily concerned with the connection defined by
% When the total space of the bundle, in this case $\text{GL}^+(n)$, is equipped with a Riemannian metric, the natural choice is to construct the horizontal sub-bundle as the orthogonal complement of the vertical sub-bundle. In the finite dimensional analogue of Otto's framework, one considers $\text{GL}^+(n)$ equipped with the (flat) Riemannian metric given in coordinates by
% \begin{align}
%     \mathcal{G}_{\factor}(\dot{\factor}_1,\dot{\factor}_2)&= \text{trace}(\dot{\factor}_1\Sigma_{\text{ref}}\dot{\factor}_2^\intercal),
% \end{align}
% where $\dot{\factor}_1,\dot{\factor}_2\in T_\factor\text{GL}^+(n)$, and $\Sigma_{\text{ref}}\in\text{Sym}^+(n)$ is the same fixed covariance appearing in the definition of the map $\pi$. The metric $\mathcal{G}$ is induced by the restriction of a flat $L_2$-metric on $\text{Diff}(\mathbb{R}^n)$ to $\text{LDiff}(\mathbb{R}^n)$ \cite{otto2001geometry,modin2016geometry}. The orthogonal complement of the vertical sub-bundle with respect to $\mathcal{G}$ is the bundle $\text{Hor}:\text{GL}^+(n)\ni \factor\mapsto \text{Hor}_{\factor}\subset T_{\factor}\text{GL}^+(n)$ defined by:
% \begin{align*}
%     \text{Hor}_{\factor}:=\{\dot{\factor}\in T_{\factor}\text{GL}^+(n)~|~ \mathcal{G}_\factor(\dot{\factor}, V) = 0, \,\forall V\in \text{Ver}_\factor\}.
% \end{align*}
% An explicit computation \cite{modin2016geometry} shows that
\begin{align}\label{eq:horizontal_distribution}
    \text{Hor}_{\factor}=\{\dot{\factor}\in T_{\factor}\text{GL}^+(n)~|~ \dot{\factor} \factor^{-1} \in\text{Sym}(n)\}.
\end{align}
Our interest in \eqref{eq:horizontal_distribution} stems from the fact it is intimately connected with the geometry of optimal mass transport. Indeed, the condition that $\dot{\factor} \factor^{-1} \in\text{Sym}(n)$ is equivalent to the requirement that $A_t$ in \eqref{eq:system} is symmetric. Consequently, the sub-bundle \eqref{eq:horizontal_distribution} encodes the infinitesmal constraint that the transition from $\factor_{t}$ to $\factor_{t+\text{d}t}$ is a Monge map. The following proposition establishes that \eqref{eq:horizontal_distribution} is, indeed, an Ehresmann connection.
\begin{prop}
    The sub-bundle \eqref{eq:horizontal_distribution} defines an Ehresmann connection on the principal bundle $\pi$.
\end{prop}
\vspace{0.2cm}
Henceforth, we refer to the connection defined by \eqref{eq:horizontal_distribution} as the OMT connection and we denote it by $\Gamma$.

\subsection{Parallel transport along Wasserstein curves}

{ We are now in a position to consider admissible curves on the space of factors via the connection that allows to transport between fibers. This entails the notion of parallel transport where curves on the base space of covariance matrices are lifted to the space of factors. We henceforth use the compact notation $\Sigma_{.},\Phi_{.}$, subscribing a dot, to indicate functions of time taking values $\Sigma_t,\Phi_t$ for $t\in\mathbb R$.}

\subsubsection*{Horizontal lifts} A $\mathcal{C}^1$ curve $\factor_{\cdot}:[0,1]\rightarrow \text{GL}^+(n)$ is said to be a \emph{horizontal curve} with respect to $\Gamma$ if $\dot{\factor}_t\in\text{Hor}_{\factor_t}$ for all $t\in[0,1]$. Accordingly, given a $\mathcal{C}^1$ curve 
\[
\Sigma_{.}\;:\;[0,1]\rightarrow\text{Sym}^+(n)\;:\;t\mapsto \Sigma_t,
\]
a \textit{horizontal lift} of $\Sigma_{.}$ with respect to $\Gamma$ starting at the initial point $\factor_{\text{in}}\in\pi^{-1}(\Sigma_0)$ is a horizontal curve 
\[
\factor_{.} \;:\;[0,1]\rightarrow\text{GL}^+(n)\;:\; t\mapsto \factor_t
\]
satisfying
\begin{align}
    \factor_0= \factor_{\text{in}}, \;\;\Sigma_t = \pi(\factor_t),
\end{align}
for $t\in[0,1]$. Existence and uniqueness of horizontal lifts is a standard result in the theory of principal bundles equipped with Ehresmann connections \cite{kobayashi1996foundations}. For the principal bundle $\pi$, we have the following {explicit construction. The proof is straightforward and hence omitted.}

\begin{prop}\label{prop:horizontal_lift}
    Let $\Sigma_{.}:[0,1]\rightarrow \text{Sym}^+(n)$ be a $\mathcal{C}^1$ curve, $A_t$ be the unique solution to the Lyapunov equation
    \begin{align}
         A_t \Sigma_t + \Sigma_t  A_t = \dot{\Sigma}_t,
    \end{align} 
    for $t\in[0,1]$,
    and let $\Phi:[0,1]\rightarrow\mathbb{R}^{n\times n}$ be the solution to
    \begin{align}\label{eq:transportPhi}
        \dot{\Phi}_t= A_t \Phi_t, \mbox{ with } \Phi_0= I.
    \end{align}
     Then, the (unique) horizontal lift of the curve $\Sigma_{.}$ with respect to $\Gamma$ starting at the initial point $\factor_{\text{in}}\in\pi^{-1}(\Sigma_0)$ exists and is of the form
    \begin{align}\label{eq:horizontal_lift}
    \hat\factor_{.} \;:\; t&\mapsto\hat\factor_t = \factor_t\factor_{\text{in}}, \mbox{ for } t\in[0,1].
    \end{align}
    \vspace{-0.1cm}
\end{prop}

\begin{remark}
    The essence of the above proposition is to relate how the particles are transported by potential forces which are designed to steer their Gaussian distribution along the curve specified by $\Sigma_t$ in the base space. In particular, the proposition implies that the connection $\Gamma$ is \emph{complete}, i.e. that every $\mathcal{C}^1$ curve in $\text{Sym}^+(n)$ has a $\mathcal{C}^1$ horizontal lift with respect to $\Gamma$ in $\text{GL}^+(n)$. It is also of interest that \eqref{eq:transportPhi} can be expressed directly in terms of $\Sigma_t,\dot\Sigma_t$, albeit not in closed form, since
\begin{align}\label{eq:sylvester_operator}
        A_t=\int_0^\infty e^{-\tau\Sigma_t}\dot\Sigma_t
    e^{-\tau\Sigma_t}d\tau =:\mathcal L_{\Sigma_t}(\dot\Sigma_t)
    \end{align}
    using standard Lyapunov theory. $\Box$
\end{remark}

\subsubsection*{Parallel transport} The horizontal lift of a $\mathcal{C}^1$ curve $\Sigma_{.} :[0,1]\rightarrow\text{Sym}^+(n)$ takes an initial point $\factor_{\text{in}}$ on the fiber $\pi^{-1}(\Sigma_0)$ to the final point $\factor_1\factor_{\text{in}}$ on $\pi^{-1}(\Sigma_1)$. By fixing the curve $\Sigma_{.}$ and varying the initial point $\factor_{\text{in}}\in\pi^{-1}(\Sigma_0)$ for the corresponding horizontal lift, we obtain a morphism
\begin{align*}
    \pi^{-1}(\Sigma_0)\rightarrow \pi^{-1}(\Sigma_1).
    %, \;\; \factor_{\text{in}}\mapsto \factor_1\factor_{\text{in}}.
\end{align*}
between fibers. This morphism, denoted by $\text{Par}[\Sigma_{.}]$, is the \textit{parallel transport map} along the curve $\Sigma_.$ and, by construction, commutes with the right action of $\text{SO}(n,\Sigma_{\text{ref}})$ \cite{kobayashi1996foundations,montgomery2002tour}. That is, if $\Theta\in \text{SO}(n,\Sigma_{\text{ref}})$, then
\begin{align}
    \text{Par}[\Sigma_{.}](\factor_{\text{in}}\Theta) = \text{Par}[\Sigma_{.}](\factor_{\text{in}})\Theta.
\end{align}
From Proposition \ref{prop:horizontal_lift}, it is clear that the parallel transport for the connection $\Gamma$ is
\begin{align}
    \text{Par}[\Sigma_{.}](\factor_{\text{in}})=\factor_1 \factor_{\text{in}}.
\end{align}
With a slight notational abuse, we identify $\text{Par}[\Sigma_{.}]$ with the $\factor_1$, i.e.\ the endpoint of the solution to \eqref{eq:transportPhi},
since $\text{Par}[\Sigma_{.}]$ amounts precisely to the corresponding linear transformation of the factors.
% Naturally, since $\text{Sym}^+(n)$ is a Riemannian manifold when equipped with the Wasserstein-Otto metric $\bar{\mathcal{G}}$ defined in \eqref{eq:riemannian_metric_normals}, it is of special interest to characterize the parallel transport map along \textit{geodesics} with respect to the metric $\bar{\mathcal{G}}$. To that end, we have the following proposition.
% \begin{prop}
%     Let $\Sigma:[0,1]\rightarrow\text{Sym}^+(n)$ be a geodesic with respect to the metric $\bar{\mathcal{G}}$. Then, the following is true.
%     \begin{enumerate}
%         \begin{subequations}
%         \item The unique horizontal lift of $\Sigma$ with respect to $\Gamma$ starting from $A_{\text{in}}\in\pi^{-1}(\Sigma_0)$ is the curve
%         \begin{align}\label{eq:horizontal_lift_of_wasserstein_geodesics}
%             A_t&= ((1-t) I + t A^\star)A_{\text{in}} , & \forall t&\in[0,1],
%         \end{align}
%         where $A^\star\in \text{GL}^+(n)$ is
%         \begin{align}
%             A^\star = A^\star&= \Sigma_{0}^{-\frac{1}{2}}\big(\Sigma_{0}^{\frac{1}{2}}\Sigma_1\Sigma_{0}^{\frac{1}{2}}\big)\Sigma_{0}^{-\frac{1}{2}}
%         \end{align}
%         \item The parallel transport map $\text{Par}[\Sigma]$ is given by
%         \begin{align}
%             \text{Par}[\Sigma](A) &= A^\star A.
%         \end{align}
%         \end{subequations}
%     \end{enumerate}
% \end{prop}
% \vspace{0.3cm}

\subsubsection*{The holonomy group} When $\Sigma_{.}$ is a \emph{loop based at $\Sigma_0$}, i.e. a closed curve in that $\Sigma_1=\Sigma_0$, then $\text{Par}[\Sigma_{.}]$ is a linear isomorphism of $\pi^{-1}(\Sigma_0)$ onto itself. Since the fibers of the bundle $\pi$ are parameterized by $\text{SO}(n,\Sigma_{\text{ref}})$, it follows that, by fixing an arbitrary element $\factor_{\text{in}}\in \pi^{-1}(\Sigma_0)$, every loop $\Sigma_{.}$ based at $\Sigma_0$ is associated with a unique $\hat{\Theta}\in\text{SO}(n,\Sigma_{\text{ref}})$ such that 
\begin{align}
    \text{Par}[\Sigma_{.}]\factor_{\text{in}} = \factor_{\text{in}}\,\hat{\Theta}.
\end{align}
If $\tilde{\factor}_{\text{in}}$ is any other element of the fiber $\pi^{-1}(\Sigma_0)$, then $\tilde{\factor}_{\text{in}}=\factor_{\text{in}}\Theta$ for some unique $\Theta\in\text{SO}(n,\Sigma_{\text{ref}})$, and so
\begin{align}
    \text{Par}[\Sigma_{.}]\tilde{\factor}_{\text{in}}=\text{Par}[\Sigma_{.}]\factor_{\text{in}}\Theta = \tilde{\factor}_{\text{in}}\Theta^{-1}\,\hat{\Theta}\Theta.
\end{align}
We observe that the matrix $\Theta^{-1}\,\hat{\Theta}\Theta$ also belongs to $\text{SO}(n,\Sigma_{\text{ref}})$ and, by definition, is conjugate to $\hat{\Theta}$. Thus, {\em the set of all parallel transport maps along loops based at $\Sigma_0$} forms a group that can be realized as a subgroup of $\text{SO}(n,\Sigma_{\text{ref}})$.  This is the \textit{holonomy group} of the connection $\Gamma$ with reference point $\Sigma_0$ and will be denoted by $\text{Hol}(\Sigma_0)$.

% \begin{remark}
The question raised at the end of the previous section, as to whether any $\Theta \in \text{SO}(n,\Sigma_{\text{ref}})$ can be realized as the holonomy of some loop, is equivalent to the question of whether $\text{Hol}(\Sigma_0)= \text{SO}(n,\Sigma_{\text{ref}})$. The classical approach to answering this question, in the context of principal bundles, is the Ambrose-Singer theorem which gives a complete characterization of the Lie algebra of $\text{Hol}(\Sigma_0)$ in terms of the \textit{curvature form} of the principal connection \cite{kobayashi1996foundations}. However, the technical machinery required to state and invoke the theorem lies beyond the scope of this manuscript. Instead, we opt to follow a control-theoretic approach to characterizing $\text{Hol}(\Sigma_0)$. This is done in the next section.

\section{Sub-Riemannian Geometry of Linear Transport}\label{sec:IV}

We recall our goals from the introduction section. Our first goal is to provide an answer to whether $\text{Hol}(\Sigma_0)= \text{SO}(n,\Sigma_{\text{ref}})$, which is equivalent to the question of whether any $\Theta \in \text{SO}(n,\Sigma_{\text{ref}})$ can be realized as the parallel transport map along a loop in $\text{Sym}^+(n)$. Our second goal is \emph{optimal synthesis}, i.e. finding the shortest curve in $\text{Sym}^+(n)$ with a specified parallel transport map. Both questions are answered in this section following a purely control-theoretic approach.

\subsection{The holonomy of linear OMT}
To achieve our first goal, we rephrase our question, and ask instead as to whether any two elements $\factor_{\text{in}}$ and $\factor_{\text{fn}}$ in $\text{GL}^+(n)$ can be joined by a horizontal curve $\Phi_.$, i.e. the horizontal lift of a curve $\Sigma_.$ in the base space. Clearly, a positive answer implies that $\text{Hol}(\Sigma)= \text{SO}(n,\Sigma_{\text{ref}})$, for every $\Sigma\in\text{Sym}^+(n)$. The following proposition provides such an answer.

\begin{prop}\label{prop:controllability}
    Let $\factor_{\text{in}},\factor_{\text{fn}}\in \text{GL}^+(n)$. Then, there exists a horizontal curve $\factor_{\cdot}:[0,1]\rightarrow \text{GL}^+(n)$ such that $\factor_0=\factor_{\text{in}}$ and $\factor_1=\factor_{\text{fn}}$.
\end{prop}
\begin{proof}
    The proposition follows directly from the Chow-Rashevskii theorem, provided that the horizontal sub-bundle is \emph{bracket-generating}, i.e.,  satisfies the H\"ormander condition \cite{agrachev2019comprehensive}. The latter assertion follows from the observation that if $A,B\in\text{Sym}(n)$, then the commutator of the two vector fields $\dot\factor^A_t:=A\factor_t$ and $\dot\factor^B_t=B\factor_t$ is
    \newcommand{\bleft}{\bm{[} \hspace*{-2pt} \bm{[}}
    \newcommand{\bright}{\bm{]} \hspace*{-2pt} \bm{]}}
    \begin{align}
        [\dot\factor^A_t, \dot\factor^B_t] = \bleft B,A\bright \factor_t = (B A - A B)\factor_t,
    \end{align}
    where $[\cdot,\cdot]$ denotes the Lie-bracket between vector fields \tcb{and $\bleft\cdot,\cdot\bright$ denotes the matrix commutator; note the apparent reversal in order, as compared to the Lie bracket.}
    Any matrix $P\in\mathbb{R}^{n\times n}$ can be written as the sum of its symmetric part, say, $S$ and its skew symmetric matrix $A$. In turn, any skew symmetric part can be written as the commutator of two symmetric matrices, i.e., 
    \begin{align}
        A= \bleft B, C\bright,
    \end{align}
    for some $B,C\in\text{Sym}(n)$.

\tcb{Before we continue with the proof of the proposition, we give a short proof of the last statement. The validity of the statement can easily be seen in the case of  
$2\times 2$ skew-symmetric matrix, since any such matrix is of the form
    \begin{align*}
        A_\lambda:=\begin{pmatrix}0 & \lambda \\-\lambda & 0 \end{pmatrix},
    \end{align*}
for some $\lambda\in\mathbb R$, while observing that $A_\lambda= \bleft \lambda E_1,E_2\bright$ for
    \begin{align*}
 E_1&=\begin{pmatrix}1 & 0 \\0 & 0 \end{pmatrix}, &
        E_2&=\begin{pmatrix}0 & 1 \\1 & 0 \end{pmatrix}.
            \end{align*}
  In general, any skew-symmetric matrix $A$ can be brought  by an orthogonal transformation $\Theta$  into a block diagonal form
  \[
  \Theta A\Theta^\top = \begin{bmatrix}
  A_{\lambda_1} &  0 & \cdots & 0 \\
  0 & A_{\lambda_2} & & 0 \\
  \vdots & & \ddots  & \vdots\\
  0 & 0 & \cdots & A_{\lambda_r}
\end{bmatrix},
  \]
  for some $\lambda_1,\ldots, \lambda_r\in \mathbb R$ when $n$ is even, or with an additional $0$ on the diagonal when $n$ is odd. (This last statement can be easily shown using the fact that any skew-symmetric matrix is normal.) Thus, in general, $A$ can be written as the commutator of the symmetric matrices
 \begin{align*}
   B=  \Theta^\top \begin{bmatrix}
  \lambda_1 E_1 &  0 & \cdots & 0 \\
  0 & \lambda_2 E_1 & & 0 \\
  \vdots & & \ddots & \vdots \\
  0 & 0 & \cdots & \lambda_r E_1
\end{bmatrix}\Theta,
 \end{align*}
 and
  \begin{align*}
   C=  \Theta^\top \begin{bmatrix}
  E_2&  0 & \cdots & 0 \\
  0 & E_2 & & 0 \\
  \vdots & & \ddots & \vdots \\
  0 & 0 & \cdots & E_2 \\
\end{bmatrix}\Theta,
 \end{align*}
 when $n$ is even, and similarly constructed matrices with an additional $0$ on the diagonal if $n$ is odd.
 }

    \tcb{Returning to the proof,} we conclude that
    \begin{align}
        \text{span}_{\mathbb{R}}\{A,\bleft B,C\bright\}_{A,B,C\in\text{Sym}(n)} = \mathbb{R}^{n\times n}.
    \end{align}
    from whence it follows that, for any $\factor\in\text{GL}^+(n)$,
    \begin{align}\label{eq:strong_bracket_condition}
        \text{span}_{\mathbb{R}}\{\dot\factor^A, [\dot\factor^B,\dot\factor^C]\}_{A,B,C\in\text{Sym}(n)} = T_{\factor}\text{GL}^+(n),
    \end{align}
    and, thus, the horizontal sub-bundle is bracket-generating.  
\end{proof}

\tcb{
As an immediate consequence of Proposition \ref{prop:controllability}, we have the following Corollary.
\begin{cor}
    $\text{Hol}(\Sigma)= \text{SO}(n,\Sigma_{\text{ref}})$ for all $\Sigma\in\text{Sym}^+(n)$.
\end{cor}
\begin{proof}
 Since Proposition \ref{prop:controllability} establishes that any two points in $\text{GL}^+(n)$ can be connected with a horizontal curve, the same applies to any two points that belong to the same fiber. These are then necessarily related by an element of the structure group, i.e. $\mathrm{SO}(n,\Sigma_{\mathrm{ref}})$. Thus, every $\Theta\in \mathrm{SO}(n,\Sigma_{\mathrm{ref}})$ can be realized as the parallel transport map along a loop in $\text{Sym}^+(n)$. The conclusion of the corollary follows.
\end{proof}}
\tcb{Proposition \ref{prop:controllability} goes beyond the statement of the Corollary since it establishes} one of the essential ingredients for defining a \emph{sub-Riemannian} metric on the principal bundle $\pi:\text{GL}^+(n)\rightarrow \text{Sym}^+(n)$; namely that the horizontal sub-bundle \eqref{eq:horizontal_distribution} is \emph{bracket-generating}. The remaining ingredient is a choice of a Riemannian metric on the base space $\text{Sym}^+(n)$ of the bundle. Naturally, the most pertinent choice in the context of OMT is the Wasserstein-Otto metric \cite{otto2001geometry,khesin2018geometric}.
% \end{remark}

\subsection{The sub-Riemannian structure of linear OMT}
Otto's metric \eqref{eq:otto}, also known as the Wasserstein-Otto metric \cite{khesin2018geometric}, specialized to Gaussian distributions, takes the form
\begin{align*}
    \langle \dot\mu_t,\dot\mu_t\rangle &=
\trace \left(A_t \Sigma_t A_t\right),
\end{align*}
%in light of \eqref{eq:potential}, 
where $A_t$ can be expressed directly as a function of $\Sigma_t,\dot\Sigma_t$, as we already noted in \eqref{eq:sylvester_operator}. Thereby, expressing the metric directly on the base space of the bundle, $\text{Sym}^+(n)$, the Wasserstein-Otto metric \cite{otto2001geometry,khesin2018geometric} takes the form
\begin{align}\label{eq:wasserstein_otto_metric}
    \bar{\mathcal{G}}_{\Sigma}(\dot{\Sigma}_1,\dot{\Sigma}_2)= \trace(\mathcal{L}_{\Sigma}(\dot{\Sigma}_1) \Sigma \mathcal{L}_{\Sigma}(\dot{\Sigma}_2)).
\end{align}
Here, again, $\Sigma\in\text{Sym}^+(n)$, $\dot{\Sigma}_1,\dot{\Sigma}_2\in T_{\Sigma}\text{Sym}^+(n)$, and $\mathcal{L}_{\Sigma}:\text{Sym}(n)\rightarrow \text{Sym}(n)$ is the operator defined in \eqref{eq:sylvester_operator}.

Since, by construction, the restriction of $\text{d}\pi$ to the horizontal sub-bundle is an isomorphism, the Wasserstein-Otto metric \eqref{eq:wasserstein_otto_metric} defines a family of inner products on the horizontal sub-bundle:
\begin{align}
    \mathcal{G}_{\factor}(\dot{\factor}_1,\dot{\factor}_2):= \bar{\mathcal{G}}_{\pi(\Phi)}(\text{d}\pi_\factor\dot{\factor}_1,\text{d}\pi_\factor\dot{\factor}_2),
\end{align}
where $\factor\in\text{GL}^+(n)$, $\dot{\factor}_1,\dot{\factor}_2\in \text{Hor}_{\factor}$. Explicit computation shows that
\begin{align}\label{eq:wasserstein_otto_subriemannian_metric}
    \mathcal{G}_{\factor}(\dot{\factor}_1,\dot{\factor}_2)= \trace(\dot{\factor}_1\Sigma_{\text{ref}}\dot{\factor}_2^\intercal).
\end{align}
Then, we have the following proposition. 

\vspace{0.2cm}
\begin{prop}\label{prop:subriemannian_structure}
    The family of inner products \eqref{eq:wasserstein_otto_subriemannian_metric}, in conjunction with the horizontal bundle $\eqref{eq:horizontal_distribution}$, define a sub-Riemannian metric on $\text{GL}^+(n)$.
\end{prop}

\vspace{0.2cm}
Henceforth, we refer to the sub-Riemannian metric defined by Proposition \ref{prop:subriemannian_structure} as the OMT sub-Riemannian metric and denote it by $\mathcal{G}$. 

\subsubsection*{Geodesics and the sub-Riemannian distance} It is a standard result that a sub-Riemannian metric induces a distance function \cite[Definition 3.30]{agrachev2019comprehensive}. So let $d_{SR}:\text{GL}^+(n)\times\text{GL}^+(n)\rightarrow\mathbb{R}_{\geq 0}$ be this sub-Riemannian distance induced by $\mathcal{G}$. A \emph{sub-Riemannian geodesic} is a horizontal curve
\[
\factor_{.}:[0,1]\rightarrow\text{GL}^+(n)
\]
such that
\begin{align}
    d_{SR}(\factor_0,\factor_1)^2=\int_0^1{\mathcal{G}_{\factor_t}(\dot{\factor}_t,\dot{\factor}_t)}\,\text{d}t.
\end{align}
\tcb{The regularity of geodesics in sub-Riemannian geometry can be quite subtle. Geodesics are typically characterized as \emph{normal} and \emph{abnormal} \cite[Definition 3.60]{agrachev2019comprehensive}. Normal geodesics are differentiable and satisfy the first-order (Pontryagin) necessary conditions for optimality.
On the hand, abnormal geodesics may or may not satisfy the first-order necessary conditions of optimality. Accordingly, there are abnormal geodesics that are also normal, and others that are not; the latter are referred to as \emph{strictly abnormal}.
As shown by Montgomery \cite[Chapter 3, Theorem 3.4]{montgomery2002tour}, the presence of strictly abnormal geodesics can be a generic phenomenon.} If present, strictly abnormal geodesics significantly complicate the local geometry of a sub-Riemannian manifold \cite[Chapter 11]{agrachev2019comprehensive}.
Fortunately, the OMT sub-Riemannian structure is special and does not admit strictly abnormal geodesics. The following proposition establishes this fact. The proof is straightforward and is sketched after the proposition. 

% There are two types of sub-Riemannian geodesics: \emph{normal} and \emph{abnormal} \cite[Definition 3.60]{agrachev2019comprehensive}. Normal geodesics are differentiable and satisfy 
% the first-order (Pontryagin) necessary conditions for optimality.
% Abnormal geodesics, on the other hand, do not satisfy the first-order necessary conditions of optimality.
% As shown by Montgomery \cite[Chapter 3, Theorem 3.4]{montgomery2002tour}, the presence of abnormal geodesics in sub-Riemannian geometry can be a generic phenomenon. If present, abnormal geodesics significantly complicate the local geometry of a sub-Riemannian manifold \cite[Chapter 11]{agrachev2019comprehensive}.
% Fortunately, the OMT sub-Riemannian structure is special and does not admit abnormal geodesics. The following proposition establishes this fact. The proof is straightforward and is sketched after the proposition. 

\vspace{0.2cm}
\begin{prop}\label{prop:no_abnormal_geodesics}
    Every nontrivial geodesic with respect to the OMT sub-Riemannian structure is normal.
\end{prop}
\vspace{0.2cm}

\begin{proof}
 From \eqref{eq:strong_bracket_condition}, we see that only first-order Lie brackets are required to generate the tangent bundle from the horizontal sub-bundle of the OMT connection. This property appears under many names in the literature. In \cite{agrachev2019comprehensive}, the horizontal sub-bundle is said to be a ``step-2" bundle when this property holds. Elsewhere, the property has been referred to as the \emph{strong bracket generating condition} \cite{strichartz1986sub}. Strichartz established in \cite[Theorem 6.1, Corollary 6.2]{strichartz1986sub} that, under the strong bracket generating condition, all sub-Riemannian geodesics are normal, see also \cite[Corollary 12.14]{agrachev2019comprehensive}. Thus, the same is true for the OMT sub-Riemannian structure. 
\end{proof}

The next proposition formulates the first-order necessary conditions for a horizontal curve to be a normal geodesic with respect to the OMT sub-Riemannian structure.
\vspace{0.2cm}
\begin{prop}\label{prop:geodesic_equations}
    Let $\factor_{.}:[0,1]\rightarrow\text{GL}^+(n)$ be a normal sub-Riemannian geodesic connecting the end points $\factor_0$ and $\factor_1$. Then, $\factor_{.}$ is smooth and there exists a unique smooth curve $\Lambda_{.}:[0,1]\rightarrow\text{GL}^+(n)$ such that
    \begin{align}
        \dot{{\factor}}_t&= A_t {\factor}_t, & \dot{{\Lambda}}_t&= -A_t({\Lambda}_t+
        2A_t{\factor}_t\Sigma_{\text{ref}}),
    \end{align}
    with $ A_t =-\frac{1}{2}\mathcal{L}_{\pi(\factor_t)}({\factor}_t {\Lambda}_t^\intercal+{\Lambda}_t{\factor}_t^\intercal)$.
\end{prop}
\vspace{0.1cm}
\begin{proof}
    The Hamiltonian associated with the geodesics of the OMT sub-Riemannian structure is given in coordinates by
    \begin{align*}
    H(\factor,\Lambda,A)&=\trace(A\factor\Sigma_{\text{ref}}\factor^\intercal A^\intercal + \Lambda^\intercal A\factor)\\
    &=\trace(A\factor\Sigma_{\text{ref}}\factor^\intercal A^\intercal + \factor \Lambda^\intercal A),
    \end{align*}
    wherein $A$ is the control input. Since $A\in\text{Sym}(n)$,
    \begin{align*}
        H({\factor},{\Lambda},A)=\trace(A{\factor}\Sigma_{\text{ref}}{\factor}^\intercal A +{\textstyle \frac{1}{2}} ({\factor} {\Lambda}^\intercal+{\Lambda}{\factor}^\intercal) A).
    \end{align*}
    The first variation of $H$ with respect to its arguments is
    \begin{align*}
        \delta H&= \trace(\delta {\Lambda}^\intercal A {\factor}) + \trace(\delta {\factor}^\intercal A ({\Lambda}+2 A{\factor}\Sigma_{\text{ref}})) \\
        &+ \trace(\delta A ({\factor}\Sigma_{\text{ref}}{\factor}^\intercal A + A{\factor}\Sigma_{\text{ref}}{\factor}^\intercal+{\textstyle \frac{1}{2}}({\factor} {\Lambda}^\intercal+{\Lambda}{\factor}^\intercal))).
    \end{align*}
    From the calculus of variations, the necessary conditions of optimality are
    \begin{align}
        \dot{{\factor}}_t&= A_t {\factor}_t, & \dot{{\Lambda}}_t&= - A_t ({\Lambda}_t+2 A_t{\factor}_t\Sigma_{\text{ref}}),
    \end{align}
    whereas Pontryagin's minimum principle necessitates that, for all $t\in[0,1]$, $A_t$ is the unique solution to the Lyapunov equation
    \begin{align}\label{eq:pmp_condition}
        {\factor}_t\Sigma_{\text{ref}}{\factor}_t^\intercal A_t + A_t{\factor}_t\Sigma_{\text{ref}}{\factor}_t^\intercal &=-{\textstyle \frac{1}{2}}({\factor}_t {\Lambda}_t^\intercal+{\Lambda}_t{\factor}_t^\intercal).
    \end{align}
    The proof is concluded by noting that ${\factor}_t\Sigma_{\text{ref}}{\factor}_t^\intercal=\pi(\factor_t)$, and that the solution to \eqref{eq:pmp_condition} is expressed via the linear operator defined in \eqref{eq:sylvester_operator}.
\end{proof}
\vspace{0.2cm}

% \subsubsection*{Invariance under the action of $\text{SO}(n,\Sigma_{\text{ref}})$}
% An important observation is that the sub-Riemannian metric $\mathcal{G}$ is invariant under the action of the structure group of the bundle, i.e. $\text{SO}(n,\Sigma_{\text{ref}})$. This fact is formalized in the following proposition. 

% \vspace{0.2cm}
% \begin{prop}
    
% \end{prop}
% \vspace{0.2cm}

% Sub-riemannian metrics that satisfy this property are said to be of the \emph{bundle type} \cite[Chapter 11]{montgomery2002tour}. They are 
\subsubsection*{Invariance with respect to $\Sigma_{\text{ref}}$} Another important feature of the OMT sub-Riemmanian structure is that it is largely independent of the choice of $\Sigma_{\text{ref}}$. The following proposition formalizes this fact.

\vspace{0.2cm}
\begin{prop}\label{prop:bundle_isomorphism}
    Let $\pi:\text{GL}^+(n)\rightarrow \text{Sym}^+(n)$ and $\tilde{\pi}:\text{GL}^+(n)\rightarrow \text{Sym}^+(n)$ be the two principal bundles defined by $\pi(\factor)= L_{\factor}(\Sigma_{\text{ref}})$ and $\tilde{\pi}(\factor)= L_{\factor}(\tilde{\Sigma}_{\text{ref}})$ for any $\Sigma_{\text{ref}},\tilde{\Sigma}_{\text{ref}} \in \text{Sym}^+(n)$. Let $\mathcal{G}$ and $\tilde{\mathcal{G}}$ be the family of inner products defining the sub-Riemannian structure on $\pi$ and $\tilde{\pi}$, respectively. Then, the bundles $\pi$ and $\tilde{\pi}$ are isometrically isomorphic.
\end{prop}
\vspace{0.1cm}
\begin{proof}
    A principal bundle isomorphism \cite[p. 53]{kobayashi1996foundations} between $\pi$ and $\tilde{\pi}$ is a diffeomorphism $\Psi:\text{GL}^+(n)\rightarrow \text{GL}^+(n)$ and a group isomorphism $\psi:\text{SO}(n,\Sigma_{\text{ref}})\rightarrow \text{SO}(n,\tilde{\Sigma}_{\text{ref}})$ such that
    \begin{align}
        \Psi(\Phi\tilde{\Theta}) = \Psi(\Phi)\psi(\Theta).
    \end{align}
    Let $\hat{\factor}\in\text{GL}^+(n)$ be any element in $\tilde{\pi}^{-1}(\Sigma_{\text{ref}})$, i.e. $\hat{\factor}\sharp \tilde{\Sigma}_{\text{ref}} = \Sigma_{\text{ref}}$, and define the maps
    \begin{align}
        \Psi(\factor) &= \factor \hat{\factor}, & \psi(\Theta) &=  \hat{\factor}^{-1}\Theta \hat{\factor}.
    \end{align}
    Then, if $\Theta\in \text{SO}(n,\Sigma_{\text{ref}})$, we have that
    \begin{align}
        \Psi(\factor\Theta) = \factor\Theta\hat{\factor} = \factor\hat{\factor}\hat{\factor}^{-1}\Theta\hat{\factor} = \Psi(\factor) \psi(\Theta).
    \end{align}
    It is clear that $\Psi$ is a diffeomorphism and that $\psi$ is a group isomorphism, which implies that the pair $\Psi$ and $\psi$ define a principal bundle isomorphism. It remains to show that $\Psi$ is also an isometry. Explicit computation gives that
    \begin{align}
        \text{d}\Psi_{\factor}\dot{\factor}= \dot{\factor} \hat{\factor}. 
    \end{align}
    Hence, we have that
    \begin{align}
        \tilde{\mathcal{G}}_{\Psi(\factor)}(\text{d}\Psi_{\factor}\dot{\factor}_1,\text{d}\Psi_{\factor}\dot{\factor}_2)  =\text{trace}(\dot{\factor}_1 \hat{\factor}\tilde{\Sigma}_{\text{ref}}\hat{\factor}^\intercal\dot{\factor}_2^\intercal).
    \end{align}
    However, by definition, $\hat{\factor}\tilde{\Sigma}_{\text{ref}}\hat{\factor}^\intercal = \Sigma_{\text{ref}} $, and so
    \begin{align}
        \mathcal{G}_{\factor}(\dot{\factor}_1, \dot{\factor}_1) = \tilde{\mathcal{G}}_{\Psi(\factor)}(\text{d}\Psi_{\factor}\dot{\factor}_1,\text{d}\Psi_{\factor}\dot{\factor}_2).
    \end{align}
    It follows that $\Psi$ is also an isometry \cite[Theorem 8.2]{strichartz1986sub}.
\end{proof}
\vspace{0.2cm}

The essence of Proposition \ref{prop:bundle_isomorphism} is that the choice of the reference covariance $\Sigma_{\text{ref}}\in\text{Sym}^+(n)$ in defining the principal bundle $\pi$ and the OMT sub-Riemannian metric is completely arbitrary. It also establishes that the holonomy groups associated to different choices of $\Sigma_{\text{ref}}$ are all conjugate to one another, i.e. that they are different representations of the same algebraic object. So, the choice $\Sigma_{\rm ref}=I$ with $\text{SO}(\Sigma_{\rm ref},n)$ being simply $\text{SO}(n)$ may be a convenient option.

\subsubsection*{The Riemannian extension} Inspection of \eqref{eq:wasserstein_otto_subriemannian_metric} reveals that this is precisely the restriction to the horizontal sub-bundle of a $\Sigma_{\text{ref}}$-weighted Frobenius inner product on each tangent space $T_{\factor}\text{GL}^+(n)$. Extending \eqref{eq:wasserstein_otto_subriemannian_metric} to the entire tangent bundle defines a (flat) Riemannian metric, henceforth referred to as the $L_2$-metric, on $\text{GL}^+(n)$. In fact, the Wasserstein-Otto metric \eqref{eq:wasserstein_otto_metric} is typically constructed by showing that the map $\pi$ is a \emph{Riemannian submersion} and that the $L_2$-metric \emph{descends} to the Wasserstein-Otto metric \eqref{eq:wasserstein_otto_metric} under $\pi$ \cite{otto2001geometry,modin2016geometry}. This observation is crucial in establishing many important regularity properties of the sub-Riemannian distance $d_{SR}$. For instance, the distance $d_{SR}$ induces the same topology as its Riemannian extension. Another consequence is that, under the strong bracket generating condition, the distance $d_{SR}$ is bounded above on compact sets by a multiple of the square root of its Riemannian extension \cite[Theorem 11.1]{strichartz1986sub}. This important fact is recapitulated in the following proposition.
\vspace{0.2cm}
\begin{prop}\label{prop:riemannian_extension_estimate}
    Let $d_{R}$ be the distance associated with the $L_2$-metric and let $K$ be a compact set in the topology induced by $d_{SR}$. Then, there exists a positive constant $c$ such that: $$d_{R}(\factor_1,\factor_2)\leq d_{SR}(\factor_1,\factor_2)\leq c d_{R}(\factor_1,\factor_2)^\frac{1}{2},$$ for all $\factor_1,\factor_2\in K$. 
\end{prop}
\vspace{0.2cm}

The significance of Proposition \ref{prop:riemannian_extension_estimate} stems from the \emph{flatness} of the $L_2$-metric\cite{otto2001geometry,modin2016geometry,khesin2018geometric}. That is, geodesics with respect to the $L_2$-metric satisfy the trivial geodesic equation
\begin{align}
    \ddot{\factor}_t&= 0,
\end{align}
which is explicitly integrable; the unique geodesic $\factor_.:[0,1]\rightarrow\text{GL}^+(n)$ between any $\factor_{\text{in}}$ and $\factor_{\text{fn}}$ in $\text{GL}^+(n)$, if it exists, is given explicitly by
\begin{align}\label{eq:L2_geodesics}
    \factor_t&= \factor_{\text{in}} + t\,(\factor_{\text{fn}}-\factor_{\text{in}}). 
\end{align}
Consequently, it can be shown that the distance $d_{R}$ has the explicit expression
\begin{align}\label{eq:L2_distance}
    d_{R}(\factor_1,\factor_2)^2=\text{trace}((\factor_1-\factor_2)\Sigma_{\text{ref}}(\factor_1-\factor_2)^\intercal).
\end{align}
\begin{figure*}
    \centering
    \begin{minipage}[t]{0.32\textwidth}
        \centering
        \includegraphics[width=1\linewidth]{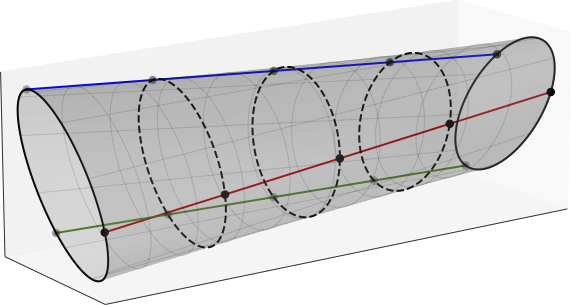}
        \caption{Trajectories of three tracer particle traversing the McCann geodesic connecting $\Sigma_{\text{in}}$ and $\Sigma_{\text{fn}}$ prescribed in \eqref{eq:prescribed_endpoints_mccann}.}
        \label{fig:mccann_curve}
    \end{minipage}
    \begin{minipage}[t]{0.32\textwidth}
        \centering
        \includegraphics[width=1\linewidth]{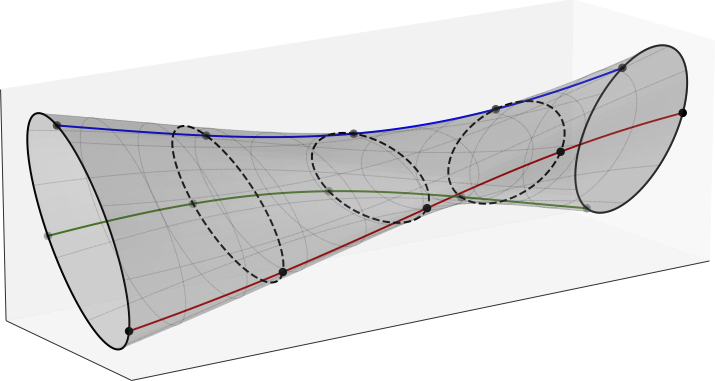} 
        \caption{Trajectories of three tracer particle traversing the isoparallel curve connecting $\Sigma_{\text{in}}$ and $\Sigma_{\text{fn}}$ prescribed in \eqref{eq:prescribed_endpoints_mccann} with $\Theta$ in \eqref{eq:prescribed_holonomy_example}. }
        \label{fig:isoparallel_curve}
    \end{minipage}
    \begin{minipage}[t]{0.32\textwidth}
        \centering
        \includegraphics[width=1\linewidth]{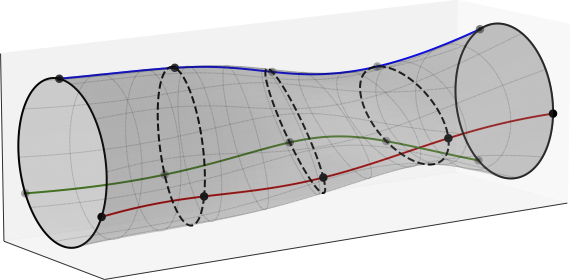}
        \caption{Trajectories of three tracer particle traversing the isoholonomic curve connecting $\Sigma_{\text{in}}=I$, $\Sigma_{\text{fn}}=I$, with $\Theta$ in \eqref{eq:prescribed_holonomy_example}.}
        \label{fig:isoholonomic_curve}
    \end{minipage}
\end{figure*}
Proposition \ref{prop:riemannian_extension_estimate}, in conjunction with \eqref{eq:L2_distance}, provide a useful upper bound for the sub-Riemannian distance $d_{SR}$. On the other hand, it is important to point out that geodesics with respect to the $L_2$-metric may fail to exist. More precisely, the curves \eqref{eq:L2_geodesics} may fail to remain in $\text{GL}^+(n)$. For example, if the dimension $n$ is even, then the curve \eqref{eq:L2_geodesics} with $\factor_{\text{in}}=I$ and $\factor_{\text{fn}}=-I$, both of which are in $\text{GL}^+(n)$, does not remain in $\text{GL}^+(n)$ for every $t\in[0,1]$. Indeed, at $t=\frac{1}{2}$, $\factor_{t}=0$! The situation is similar when $n$ is odd. This behavior is primarily due to the incompatibility of the $L_2$-metric with the group structure on $\text{GL}^+(n)$. 

\subsection{The Isoparallel Mass Transport Problem}

Having defined the OMT sub-Riemannian structure, we are in a position to pursue our second goal, the design of a control protocol that realizes a specified parallel transport map. To this end, we formulate the following problem.

\vspace{0.2cm}
\begin{prob}[Isoparallel Mass Transport (IMT)]\label{prob:isoparalell_problem}
    Given $\Sigma_{\text{in}}$, $\Sigma_{\text{fn}}$ in $\text{Sym}^+(n)$, and $\factor_{\text{des}}$ in $\text{GL}^+(n)$ with $\factor_{\text{des}}\sharp \Sigma_{\text{in}}=\Sigma_{\text{fn}}$, find a curve of minimal length $\Sigma^\star_{\cdot}:[0,1]\rightarrow\text{Sym}^+(n)$ that satisfies $\Sigma^\star_0= \Sigma_{\text{in}},$ $ \Sigma^\star_1=\Sigma_{\text{fn}},$ and $\text{Par}[\Sigma_{.}^\star] =  \Phi_{\text{des}}$.
\end{prob}
\vspace{0.2cm}

Evidently, the essence of Problem \ref{prob:isoparalell_problem} is to effect transport that links two distributions on the base space $\mathcal N(\mathbb R^n)$, corresponding to $\Sigma_{\rm in}$ and $\Sigma_{\rm fin}$, respectively, ensuring a specified terminal positioning of the particles while minimizing the quadratic cost that is associated with the length that is being traversed under the control of potential forces. \tcb{The constraint of the terminal positioning of the particles marks a departure from the classical optimal mass transport problem where only the two terminal distributions are specified.}

The IMT problem is an instance of Montgomery's \emph{isoparallel problem} \cite{montgomery1990isoholonomic} expressed in the context of OMT. It is equivalent to the problem of finding a geodesic with respect to the OMT sub-Riemannian structure. 
% Unlike the general case of the latter however, sub-Riemannian geodesics associated to solutions of the isoparallel problem are never abnormal. For a proof of this fact, we refer the reader to \cite[Theorem 1]{montgomery1990isoholonomic}. In fact, Proposition \ref{prop:no_abnormal_geodesics} already rules out the existence of nontrivial abnormal geodesics for the OMT sub-Riemannian structure, which implies that solutions to the IMT problem, when they exist, are necessarily projections of normal sub-Riemannian geodesics.

\subsubsection*{First-order necessary conditions} The following proposition is a statement of the necessary conditions satisfied by a solution of the IMT problem.
\vspace{0.2cm}
\begin{prop}\label{prop:imt_necessary_conditions}
    Assume that $\Sigma_{\cdot}^\star:[0,1]\rightarrow \text{Sym}^+(n)$ is a solution to the IMT problem. Then, there exists $\Omega\in\mathbb{R}^{n\times n}$, skew symmetric, such that $\Sigma_{\cdot}^\star$ satisfies
    % \begin{subequations}
        \begin{align}\label{eq:IMT_necessary_conditions}
            % &\dot{\Sigma}^\star_t= \mathcal{L}_{\Sigma_t^\star}(\Pi_t) \Sigma_t^\star+\Sigma_t^\star \mathcal{L}_{\Sigma_t^\star}(\Pi_t),\\ 
            % &\begin{aligned}
            %     \dot{\Pi}_t&= (\Sigma^\star_t \mathcal{L}_{\Sigma_t^\star}(\Pi_t)+\Omega)\mathcal{L}_{\Sigma_t^\star}(\Pi_t)\\
            %     &\,+\mathcal{L}_{\Sigma_t^\star}(\Pi_t)(\mathcal{L}_{\Sigma_t^\star}(\Pi_t)\Sigma^\star_t-\Omega).
            % \end{aligned}
            \ddot{\Sigma}^\star_t-(\Sigma^\star_t \mathcal{L}_{\Sigma_t^\star}(\dot{\Sigma}^\star_t)^2+\mathcal{L}_{\Sigma_t^\star}(\dot{\Sigma}^\star_t)^2\Sigma^\star_t)&= \xi(\Sigma_t^\star,\dot{\Sigma}_t^\star),
        \end{align}
        for all $t\in[0,1]$, with $\xi(\Sigma,\dot{\Sigma})=\Omega \mathcal{L}_{\Sigma}(\dot{\Sigma})-\mathcal{L}_{\Sigma}(\dot{\Sigma})\Omega$. 
    % \end{subequations}
\end{prop}
\vspace{0.2cm}
\begin{proof}
    % As noted earlier, without loss of generality, we take the reference covariance to be $\Sigma_{\text{ref}}=\Sigma_{\text{in}}$, which implies that $I\in\pi^{-1}(\Sigma_{\text{in}})$. 
    Fixing $\Phi_{\text{in}}\in\pi^{-1}(\Sigma_{\text{in}})$, the IMT problem is equivalent to the problem of finding the normal sub-Riemannian geodesic between $\Phi_{\text{in}}$ and $\factor_{\text{des}}\Phi_{\text{in}}$. From Proposition \ref{prop:geodesic_equations}, there must exist a unique smooth curve $\Lambda_.:[0,1]\rightarrow\text{GL}^+(n)$ such that the sought-after geodesic satisfies
   \begin{align*}
        \dot{{\factor}}_t&= A_t {\factor}_t, & \dot{{\Lambda}}_t&= -A_t({\Lambda}_t+2
        A_t{\factor}_t\Sigma_{\text{ref}}),
    \end{align*}
    with $A_t = -\frac{1}{2} \mathcal{L}_{\pi(\factor_t)}({\factor}_t {\Lambda}_t^\intercal+{\Lambda}_t{\factor}_t^\intercal)$. Projecting $\dot{\factor}_t$ onto the tangent space of $\text{Sym}^+(n)$, we obtain
    \begin{align*}
        \dot{\Sigma}^\star_t=\text{d}\pi_{\factor_t}\dot{\factor}_t = A_t \Sigma^\star_t+\Sigma^\star_t A_t = -\frac{1}{2} ({\factor}_t {\Lambda}_t^\intercal+{\Lambda}_t{\factor}_t^\intercal).
    \end{align*}
    Let $\Pi_t$ and $\Omega_t$ be given by
    \begin{align*}
        \Pi_t&= -\frac{1}{2}(\factor_t\Lambda_t^\intercal+\Lambda_t\factor_t^\intercal), & \Omega_t&= \frac{1}{2}(\factor_t\Lambda_t^\intercal-\Lambda_t\factor_t^\intercal).
    \end{align*}
    Then, an explicit calculation shows that
    \begin{align*}
        \dot{\Pi}_t&= (\Sigma^\star_t A_t+\Omega_t)A_t+A_t(A_t\Sigma^\star_t +\Omega_t^\intercal), & \dot{\Omega}_t&= 0.
        % \ddot{\Sigma}^\star_t&= \dot{\Pi}_t, & \dot{\Omega}_t&= 0.
    \end{align*}
    In particular, we have that
    \begin{align*}
        \Omega_t&= \Omega_0=\Omega = \frac{1}{2}(\Phi_0\Lambda_0^\intercal-\Lambda_0\Phi_0^\intercal).
    \end{align*}
    The proof is concluded by observing that $\Pi_t=\dot{\Sigma}_t^\star$. 
\end{proof}
\vspace{0.2cm}
The usefulness of Proposition \ref{prop:imt_necessary_conditions}, compared to Proposition \ref{prop:geodesic_equations}, is that $\Sigma_{\text{ref}}$ does not appear in \eqref{eq:IMT_necessary_conditions}.
% We remark that \eqref{eq:IMT_necessary_conditions} is \emph{Wong's equation} for the OMT connection \cite[Chapter 12]{montgomery2002tour}. In particular, the left hand side of \eqref{eq:IMT_necessary_conditions} is the coordinate expression for the Levi-Civita connection of the Wasserstein-Otto metric \eqref{eq:wasserstein_otto_metric}. 
\vspace{0.2cm}

\subsubsection*{Existence of solutions to the IMT problem} Sub-Riemannian geodesics between points that are sufficiently close, with respect to the sub-Riemannian distance $d_{SR}$, are guaranteed to exist \cite[Theorem 1.18]{montgomery2002tour}. The standard sufficient condition, see \cite[Proposition 1]{montgomery1990isoholonomic}, to guarantee global existence of sub-Riemannian geodesics, and consequently, global existence of solutions to the isoparallel problem, include the condition that the Riemannian structure on the base space of the principal bundle is complete, i.e. that $\text{Sym}^+(n)$ equipped with the Wasserstein-Otto metric \eqref{eq:wasserstein_otto_metric}, abbreviated as $(\text{Sym}^+(n),\bar{\mathcal{G}})$ for short, is a complete Riemannian manifold. Unfortunately, it is well-known that $(\text{Sym}^+(n),\bar{\mathcal{G}})$ is not complete. As shown in Remark \ref{rem:remark1}, nontrivial geodesics in $(\text{Sym}^+(n),\bar{\mathcal{G}})$ may be extended in one but never in both directions. This alone is sufficient to prove that the sub-Riemannian structure of OMT is not geodesically complete. Indeed, we have the following proposition which follows from the fact that the OMT sub-Riemannian structure is of the \emph{bundle-type} \cite[Chapter 11]{montgomery2002tour}.
\begin{prop}\label{prop:mccann_geodesics_are_subriemannian}
    Any horizontal lift of a McCann geodesic is a sub-Riemannian geodesic.
\end{prop}

\vspace{0.05cm}
A geodesic in $(\text{Sym}^+(n),\bar{\mathcal{G}})$ that cannot be extended must approach the boundary of $\text{Sym}^+(n)$ which consists of degenerate symmetric positive semi-definite matrices. The horizontal lift of such a curve will necessarily approach a degenerate $\Phi\in\mathbb{R}^{n\times n}\backslash\text{GL}^+(n)$, and so cannot be extended indefinitely in $\text{GL}^+(n)$. A potential remedy might be to \emph{complete} $(\text{Sym}^+(n),\bar{\mathcal{G}})$, i.e. include all its limit points in the topology induced by the metric $\bar{\mathcal{G}}$. Unfortunately, this cannot be done smoothly. Namely, by completing $(\text{Sym}^+(n),\bar{\mathcal{G}})$, it ceases to be a smooth manifold altogether. Instead, it becomes a \emph{stratified space}, though it remains a \emph{length space} \cite{takatsu2011wasserstein}. Alternatively, we may complete the OMT sub-Riemannian structure itself as a metric space. This approach is faced with a major obstacle: the horizontal sub-bundle, as defined, is not bracket-generating at the boundary of $\text{GL}^+(n)$. A third alternative is to ``regularize" the OMT sub-Riemannian metric in the spirit of the Schr\"odinger bridge approach to regularized OMT \cite{leonard2012schrodinger,chen2016relation}. The aforementioned potential remedies and their implications for the OMT sub-Riemannian structure will be explored in future work.

\subsection{Examples: sub-Riemannian geodesics}
 We conclude the current section with representative numerical examples that illustrate our results. We fix $n=2$, and $\Sigma_{\text{ref}}= I$, which implies that the isotropy group $\text{SO}(n,\Sigma_{\text{ref}})$ is simply the group of planar rotations $\text{SO}(2)$. It also implies that the fibers of the bundle are
\begin{align*}
    \pi^{-1}(\Sigma) &= \{\Sigma^{\frac{1}{2}}\Theta~|~\Theta\in\text{SO}(2)\}, & &\text{for all $\Sigma\in\text{Sym}^+(n)$.}
\end{align*} 

\subsubsection*{McCann geodesics} {\color{black}
If $\Phi_{\text{des}}=\Phi^\star$ is the matrix associated with the Monge map, i.e. the matrix given by 
\begin{align}
    \Phi^\star= \Sigma_{\text{in}}^{-\frac{1}{2}}\left(\Sigma_{\text{in}}^{\frac{1}{2}}\Sigma_{\text{fn}}^{\vphantom{\frac{1}{2}}}\Sigma_{\text{in}}^{\frac{1}{2}}\right)^{\frac{1}{2}}\Sigma_{\text{in}}^{-\frac{1}{2}},
\end{align}
then the associated sub-Riemannian geodesics are nothing but horizontal lifts of the McCann geodesic connecting the two points $\Sigma_{\text{in}}$ and $\Sigma_{\text{fn}}$ in the base space $\text{Sym}^+(n)$. A sample trajectory is shown in Fig. \ref{fig:mccann_curve} wherein the endpoints are
\begin{align}\label{eq:prescribed_endpoints_mccann}
    \Sigma_{\text{in}}&= \begin{bmatrix}3 & 2 \\ 2 & 3\end{bmatrix}, & \Sigma_{\text{fn}}&= \begin{bmatrix}\phantom{-}3 & -2 \\ -2 & \phantom{-}3\end{bmatrix}.
\end{align}
The ellipses in Fig. \ref{fig:mccann_curve}, as well as in Fig. \ref{fig:isoparallel_curve} and Fig. \ref{fig:isoholonomic_curve}, mark the evolution of the level set of the gaussian distribution to which the indicated tracer particles belong.
As can be observed, the tracer particles travel along straight lines (geodesics of $\mathbb{R}^n$) which is a characteristic feature of McCann geodesics.
}
\vspace{0.2cm}
\subsubsection*{Isoparallel geodesics}

If $\Phi_{\text{des}} = \Phi^\star \Theta$ for some $\Theta\in\text{SO}(2)$, then the associated sub-Riemannian geodesics are no longer lifts of the McCann geodesic. For $\Theta\in\text{SO}(2)$ close enough to $I$, one can show that the geodesics exist. A sample trajectory is shown in Fig. \ref{fig:isoparallel_curve} wherein the end points $\Sigma_{\text{in}}$ and $\Sigma_{\text{fn}}$ are those in \eqref{eq:prescribed_endpoints_mccann} and $\Theta$ is
\begin{align}\label{eq:prescribed_holonomy_example}
 \Theta&=\exp(\omega), & \omega &= \begin{bmatrix}
        \phantom{-}0 & \frac{1}{2}\\ -\frac{1}{2} & 0
    \end{bmatrix},
\end{align}
which corresponds to a planar rotation of $\approx 28.6^\circ$.
\vspace{0.2cm}

\subsubsection*{Isoholonomic geodesics}
If $\Sigma_{\text{fn}}=\Sigma_{\text{in}}=I$, then $\Phi_{\text{des}}=\Theta\in\text{SO}(2)$, and a sub-Riemannian geodesic realizing the parallel transport map $\Phi_{\text{des}}$ will necessarily connect two points of the same fiber $\pi^{-1}(I)=\text{SO}(2)$. In fact, there's an entire 1-parameter family of sub-Riemannian geodesics that realize $\Phi_{\text{des}}$ for every $\Phi_{\text{des}}\in\text{SO}(2)$ sufficiently close to $I$! Such family has a simple characterization in terms of any of its elements. Namely, if $\Sigma_{.}$ is a solution of the IMT problem with the above data, then $\tilde{\Theta}\sharp \Sigma_{.}$ is also a solution for every $\tilde{\Theta}\in\text{SO}(2)$. This phenomenon is ubiquitous in sub-Riemannian geometry \cite[Section 1.8]{montgomery2002tour}. A sample trajectory, with $\Theta$ as prescribed in \eqref{eq:prescribed_holonomy_example}, is shown in Fig. \ref{fig:isoholonomic_curve}. As can be observed, the tracer particles have rotated an angle of $\approx 28.6^\circ$ with respect to the origin after traversing the curve once. 

\section{Application: Transportation Cycles} \label{sec:Applications}

Modern day control applications of multiagent and particle systems have brought to the fore a slew of new type of performance objectives, such as maintaining formation, avoiding congestion, and so on, that require sophisticated path planning strategies. A model example for the theory that we put forth is the rudimentary paradigm of  particles being transported along cycles, where it is of essence that they retain their initial configuration upon returning to their starting configuration.
It is envisioned that elements of this theory will soon be tailored to more general dynamics, but for now we consider that the collection of particles are actuated by potential forces.

Our problem requires that a collection of particles
cycle through a sequence $\{\bm{\mu}_1,\bm{\mu}_2,\dots,\bm{\mu}_N\}$ of specified distributions. The control cost is quadratic, and is quantified by the length of the curve in the Wasserstein space traversed by the ensemble in their flight between target distributions. 

In contrast to the formalism of uncertainty control where the focus is on regulating distributions, herein the control protocol needs to guide individual particles along closed orbits. To this end, it is of essence that their position is ``registered'' relative to one another and relative to their common distribution. In other words, before designing a control law, one needs to know in advance how the particles are to be distributed when they arrive at the $k$th point of the cycle.

% \begin{figure}[tb]
%     \centering
%     \includegraphics[width=1\linewidth]{sub-RiemannianGeodesicIsoholonomic}
%     \caption{Trajectories of three tracer particle traversing the isoholonomic curve connecting $\Sigma_{\text{in}}=I$, $\Sigma_{\text{fn}}=I$, with the holonomy prescribed in \eqref{eq:prescribed_holonomy_example}. \tcb{The ellipses mark the evolution of the level set of the gaussian distribution to which the indicated tracer particles belong.}}
%     \label{fig:isoholonomic_curve}
% \end{figure}

The need to register (positional) data is not unique to ensemble control.
It is a problem that is encountered in many fields. Broadly speaking, registration refers to ways that one can correspond data between different collections of such. In other words, registration can be seen as the process to generate labels for the elements of a collection, so as to establish correspondence between elements of several different collections. Herein, labels signify positional correspondence.
In another instance, imagine pictures of the same object taken from slightly different perspectives, where we wish to establish automatically a correspondence between pixels in each. Over the past twenty years, OMT has been extensively used for these as well as other related problems of matching distributions \cite{haker_zhu,haker2003monge,peyre2019computational}.

% An important problem across many fields is ``registration.'' This refers to ways on how to correspond data between different collections of such. In other words, registration can be seen as the process to generate labels for the elements of a collection, so as to establish correspondence between elements of several different, yet related, collections.

% Imagine for instance pictures of the same object taken from slightly different perspectives, where we wish to establish automatically a correspondence between pixels in each. Over the past twenty years, OMT has been extensively used for that as well as other related problems of matching distributions \cite{haker_zhu,haker2003monge,peyre2019computational}. 

A natural approach for registration is to identify a reference distribution. Then, for instance, one may use optimal transportation to establish pairwise matching at the level of particles. More concretely, given the collection $\{\bm{\mu}_1,\bm{\mu}_2,\dots,\bm{\mu}_N\}$, and a reference distribution $\mu_{\text{ref}}$, ``registration" between $\bm{\mu}_i$ and $\bm{\mu}_j$ can be effected by
\begin{align}\label{eq:ref}
    \varphi_{i\mapsto j} &:= \varphi^\star_{\text{ref}\mapsto j}\circ\varphi^\star_{i \mapsto \text{ref}},
\end{align}
where $\varphi^\star_{\text{ref}\mapsto i}$ is the Monge map taking $\mu_{\text{ref}}$ to $\bm{\mu}_{i}$ and $\varphi^\star_{i\mapsto \text{ref}}$ is its inverse. Evidently, the holonomy of the loop
\[\bm{\mu}_i\xrightarrow{\varphi^\star_{i \mapsto \text{ref}}}\mu_{\rm ref}\xrightarrow{\varphi^\star_{\text{ref}\mapsto j}}\bm{\mu}_j\xrightarrow{\varphi^\star_{j\mapsto i}}\bm{\mu}_i,
\]
may not be trivial, and therefore, the choice $\varphi^\star_{j\mapsto i}$ is generally not an option. 
Instead, the pairwise matching in
\eqref{eq:ref} ensures that the holonomy of triangle 
\[\bm{\mu}_i\xrightarrow{\varphi^\star_{i \mapsto \text{ref}}}\mu_{\rm ref}\xrightarrow{\varphi^\star_{\text{ref}\mapsto j}}\bm{\mu}_j\xrightarrow{\varphi_{j\mapsto i}}\bm{\mu}_i
\]
is trivial.
Hence, the
interpolating curve between 
$\bm{\mu}_i$, $\bm{\mu}_j$ traversed by particles via an optimal control law, would constitute the solution to the IMT problem, i.e., it will be the curve of minimal sub-Riemannian length that realizes \eqref{eq:ref}.  Both the specific curve and the resulting distance depend on the choice of $\mu_{\rm ref}$. \\

% In general, the problem to select an optimal $\mu_{\rm ref}$ so that the ``Wasserstein perimeter" of the closed orbit that sequentially visits  $\{\bm{\mu}_1,\bm{\mu}_2,\dots,\bm{\mu}_N\}$ is minimal, is open, in that it can only be approached numerically. 

\subsubsection*{Wasserstein Polygons}

We consider $N$ Gaussian distributions with corresponding covariances $\{\mathbf \Sigma_1,\mathbf\Sigma_2,\dots,\mathbf\Sigma_N\}\subset\text{Sym}^+(n)$. In a consistent manner, $\Sigma_{\text{ref}}$ signifies the covariance of a reference distribution.
We consider separately the ``triangles'' formed with ``vertices''  $\mathbf{\Sigma}_i$, $\mathbf{\Sigma}_j$ and $\Sigma_{\rm ref}$, for $(i,j)\in\{(1,2),(2,3),\dots,(N,1)\}$.

Denoting by $\factor^\star_{\text{ref}\rightarrow i}$ the Monge map linking ${\Sigma}_{\text{ref}}$ to $\bm{\Sigma}_i$ and by $\factor^\star_{i\rightarrow \text{ref}}$ its inverse, we associate to every pair ($\mathbf{\Sigma}_i,\mathbf{\Sigma}_j)$ a family of Wasserstein triangles of the form
\begin{align*}
    \triangle^{i,j}_t&:=\begin{cases}
        (I + 3t (\factor^\star_{\text{ref}\rightarrow i}-I)){\sharp}{\Sigma}_{\text{ref}}, & t\in\left[0,\frac{1}{3}\right),\\
        \Sigma^{i\mapsto j}_\tau, \text{ with $\tau=3t-1$,} & t\in\left[\frac{1}{3},\frac{2}{3}\right),\\
        (I + (3t-2)(\factor^\star_{j\rightarrow \text{ref}}-I)){\sharp}\mathbf{\Sigma}_j, & t\in\left[\frac{2}{3},1\right].\\
    \end{cases}
\end{align*}
Note that two of the ``edges'' are formed by McCann geodesics, whereas the third ``edge'' of the triangle is $\Sigma^{i\mapsto j}_{.}:[0,1]\rightarrow \text{Sym}^+(n)$, and this is considered to be any $\mathcal{C}^1$ curve satisfying
\begin{align}\label{eq:vertices_condition}
    \Sigma^{i\mapsto j}_0&= \mathbf{\Sigma}_i, \;\; \Sigma^{i\mapsto j}_1= \mathbf{\Sigma}_j.
\end{align}
In particular, $\Sigma^{i\mapsto j}_{.}$ is the curve traversed by the distribution of the particles in their flight from $\mathbf{\Sigma}_i$ to $\mathbf{\Sigma}_j$. From \eqref{eq:vertices_condition}, it is clear that $\triangle^{i,j}_.$ is a closed curve in $\text{Sym}^+(n)$. As can be verified, the holonomy of $\triangle^{i,j}_{.}$ is
\begin{align*}
    \text{Par}[\triangle^{i,j}_.] = \factor^\star_{j\rightarrow \text{ref}}\text{Par}[\Sigma^{i\mapsto j}_.]\factor^\star_{\text{ref}\rightarrow i},
\end{align*}
where $\text{Par}[\Sigma^{i\mapsto j}_.]$ is the parallel transport map along $\Sigma^{ij}_{.}$. If no mixing is to take place, i.e. $\text{Par}[\triangle^{i,j}_.]=I$, then we need that
\begin{align}\label{eq:isoholonomic_triangle_condition}
    \text{Par}[\Sigma^{i\mapsto j}_.] = \factor^\star_{\text{ref}\rightarrow j}\factor^\star_{i\rightarrow \text{ref}}.
\end{align}
For any collection of curves $\Sigma^{i\mapsto j}_.$ with $(i,j)\in\{(1,2),(2,3),\dots,(N,1)\}$, we may construct the closed curve
\begin{align}\label{eq:wasserstein_polygon_boundary}
    \Sigma_t&:=\begin{cases}
        \Sigma^{1\mapsto 2}_\tau, \text{ with $\tau= Nt$,} & t\in\left[0,\frac{1}{N}\right],\\
        \Sigma^{2\mapsto 3}_\tau, \text{ with $\tau= Nt-1$,} & t\in\left[\frac{1}{N},\frac{2}{N}\right],\\
        ~~\vdots ~~ & \hfill \vdots \hfill \\
        \Sigma^{N\mapsto 1}_\tau, \text{ with $\tau= N(t-1)+1$,} & t\in\left[\frac{N-1}{N},1\right].\\
    \end{cases}
\end{align}
which sequentially visits the vertices of the Wasserstein polygon in the order 
$$1\rightarrow 2\rightarrow \dots \rightarrow N\rightarrow 1.$$
We then have the following proposition.

\begin{figure*}
    \centering
    \begin{minipage}[t]{0.3\textwidth}
        \centering
        \includegraphics[width=1\linewidth]{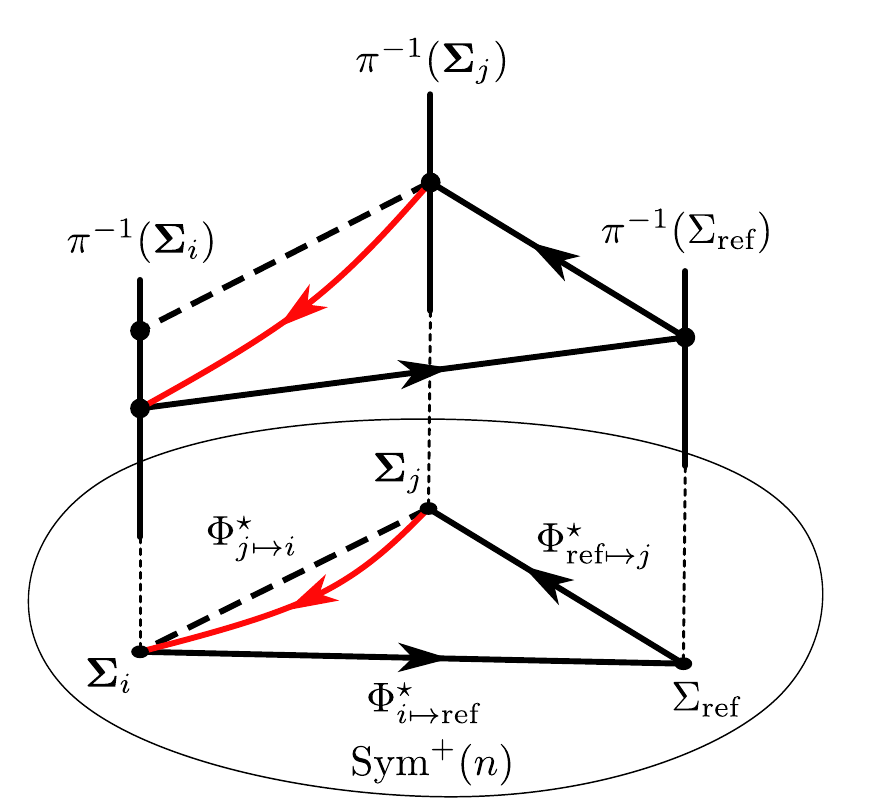}
        \caption{A sketch of Problem \ref{prob:isoholonomic_triangle}. The red curve represents a potential solution and a horizontal lift thereof. The dashed line represents a McCann geodesic between $\mathbf{\Sigma}_i$ and $\mathbf{\Sigma}_j$, and its horizontal lift. The solid black lines also indicate McCann geodesics to and from the reference measure. }
        \label{fig:isoholonomic_triangle}
    \end{minipage}
    \hfill
    \begin{minipage}[t]{0.25\textwidth}
        \centering
        \includegraphics[width=1\linewidth]{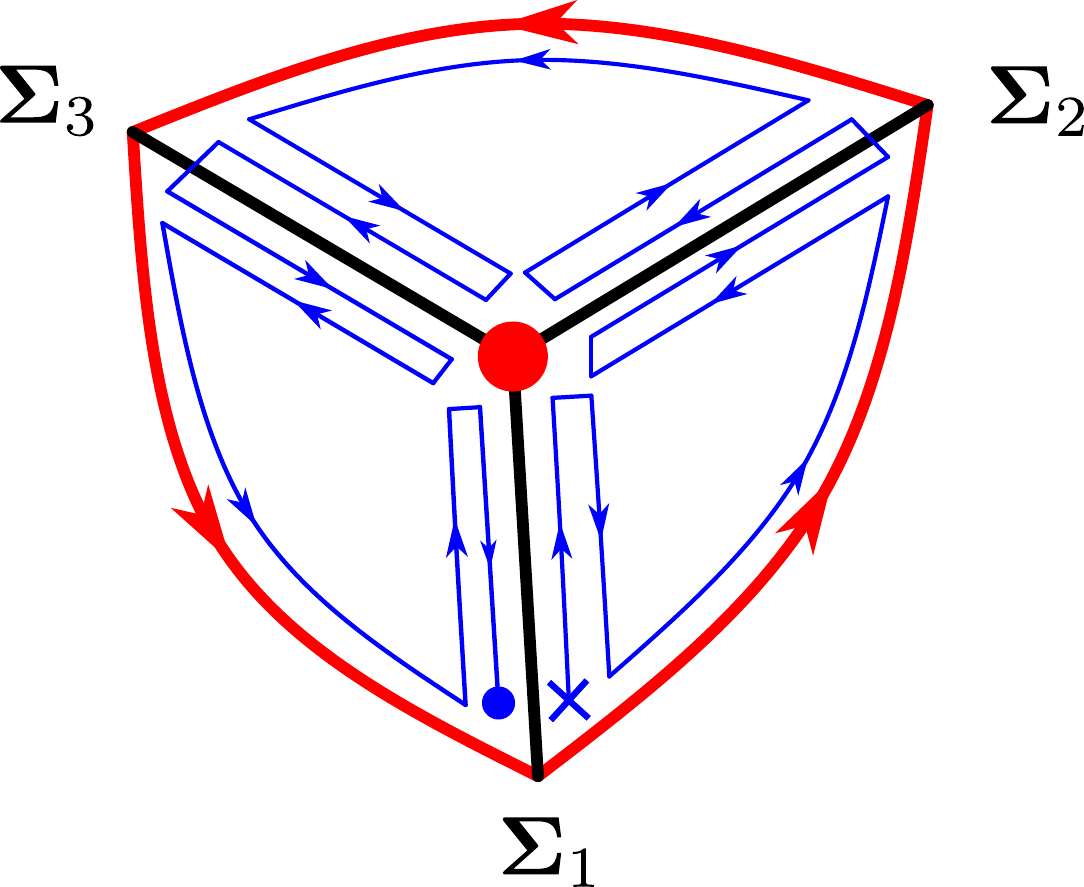}
        \caption{A geometric proof for Proposition \ref{prop:polygon_holonomy}. The red curve is $\Sigma_.$. The red circle represents $\Sigma_{\text{ref}}$. The black lines represent McCann geodesics. The blue curve has the same holonomy as $\Sigma_.$. Note that the blue curve is exaggerated for illustration. }
        \label{fig:triangle_holonomy}
    \end{minipage}
    \hfill
    \begin{minipage}[t]{.40\textwidth}
        \centering
        \includegraphics[width=1\linewidth]{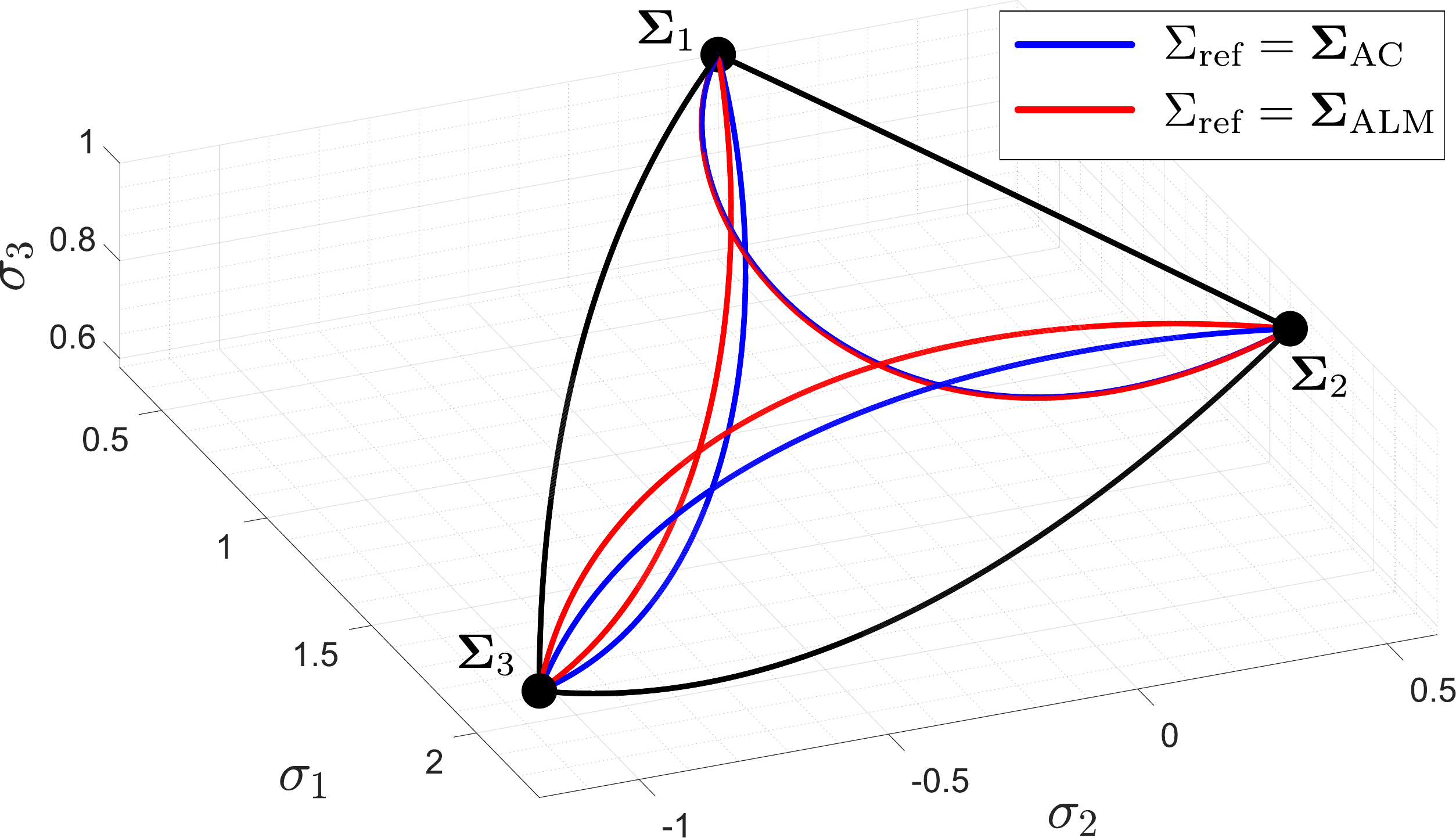}
        \caption{The isoholonomic Wasserstein triangle. The blue curve corresponds to the choice of the Agueh-Carlier Wasserstein barycenter as the reference covariance. The red curve corresponds to the choice of the Ando-Li-Mathias geometric mean as the reference covariance. The black curve represents the ``McCann triangle".}
        \label{fig:isoholonomic_triangle_combined}
    \end{minipage} 
    \hfill
\end{figure*}

\vspace{0.2cm}
\begin{prop}\label{prop:polygon_holonomy}
    If $\Sigma^{i\mapsto j}_.$ satisfies \eqref{eq:isoholonomic_triangle_condition} for every $(i,j)\in\{(1,2),(2,3),\dots,(N,1)\}$, then $\text{Par}[\Sigma_.] = I$.
\end{prop}
\vspace{0.2cm}
\begin{proof}
    By construction of the curve $\Sigma_.$, we have
    \begin{align*}
        \text{Par}[\Sigma_.]&= \text{Par}[\Sigma^{N\mapsto 1}_.]\prod_{i=1}^{N-1} \text{Par}[\Sigma^{i\mapsto i+1}_.].
    \end{align*}
    Additionally, we recall that
    \begin{align*}
        \text{Par}[\triangle^{i,j}_.] = \factor^\star_{j\rightarrow \text{ref}}\text{Par}[\Sigma^{i\mapsto j}_.]\factor^\star_{\text{ref}\rightarrow i},
    \end{align*}
    Hence, we obtain that
    \begin{align*}
        \text{Par}[\Sigma_.]&= \text{Par}[\Sigma^{N\mapsto 1}_.]\prod_{i=1}^{N-1} \factor^\star_{\text{ref}\rightarrow i+1}\text{Par}[\triangle^{i,i+1}_.]\factor^\star_{i\rightarrow \text{ref}}.
    \end{align*}
    However, by assumption, $\text{Par}[\triangle^{i,j}_.]=I$ for every pair $(i,j)\in\{(1,2),(2,3),\dots,(N,1)\}$, and so
    \begin{align*}
        \text{Par}[\Sigma_.]&= \factor^\star_{\text{ref}\rightarrow 1}\factor^\star_{N\rightarrow \text{ref}}\prod_{i=1}^{N-1} \factor^\star_{\text{ref}\rightarrow i+1}\factor^\star_{i\rightarrow \text{ref}}.
    \end{align*}
    The proposition follows by expanding the product and sequentially applying the identity $\factor^\star_{\text{ref}\rightarrow i}\factor^\star_{i\rightarrow \text{ref}} = I.$ An alternative (geometric) proof for the simple case $N=3$ is sketched in Figure \ref{fig:triangle_holonomy}. The case of arbitrary $N$ is similar.
\end{proof}

\vspace{0.2cm}

This naturally leads us to consider:
\vspace{0.2cm}
\begin{prob}\label{prob:isoholonomic_triangle}
    Find the shortest $\Sigma^{i\mapsto j}_.$ satisfying \eqref{eq:vertices_condition} and \eqref{eq:isoholonomic_triangle_condition}.   
\end{prob}
\vspace{0.2cm}
    
Clearly, Problem \ref{prob:isoholonomic_triangle} is equivalent to an IMT problem wherein
\begin{align}
    \Sigma_{\text{in}}&= \mathbf{\Sigma}_i, & \Sigma_{\text{fn}}&=\mathbf{\Sigma}_j, & \factor_{\text{des}}&= \factor^\star_{\text{ref}\rightarrow j}\factor^\star_{i\rightarrow \text{ref}}.
\end{align}
A schematic that compares the solution to the above problem with the McCann ``edge'' $(\mathbf{\Sigma}_i,\mathbf{\Sigma}_j)$ is shown in Figure \ref{fig:isoholonomic_triangle}. 

Returning to the Wasserstein polygon with vertices
$\{\mathbf \Sigma_1,\mathbf\Sigma_2,\dots,\mathbf\Sigma_N\}$, it is evident that the choice of $\Sigma_{\rm ref}$ dictates the sub-Riemannian length of the edges, and hence the sub-Riemannian ``perimeter" of the polygon, which represents the control cost in steering the ensemble along the cycle, as specified. Below we explore and compare two different options for choosing $\Sigma_{\text{ref}}$. Optimization of the control cost over $\Sigma_{\rm ref}$ will be tackled in future work. 
\tcb{\begin{remark}
    In principle, it is possible to seek a horizontal curve of minimal length that visits the given distributions, in the specified order, which leads to a sequence of optimal control problems that are nonlinearly coupled through their boundaries. This problem is numerically demanding, and the approach of registering the given distributions with respect to a reference decouples the problem into simpler ones. $\Box$
\end{remark}}

\subsubsection*{The Agueh-Carlier barycenter}
A popular choice for ``averaging" a collection of distributions $\{\bm{\mu}_1,\bm{\mu}_2,\dots,\bm{\mu}_N\}$ is the Agueh-Carlier Wasserstein barycenter \cite{agueh2011barycenters}, defined through the variational characterization
\[
    \bm{\mu}_{\text{AC}}=\arg\min_{\mu\in\mathcal{P}_2(\mathbb{R}^2)} \sum_{i=1}^N \mathcal W_2^2(\mu,\bm{\mu}_i)^2.
\]
For a collection of Gaussian measures with corresponding covariances $\{\bm{\Sigma}_1,\bm{\Sigma}_2,\dots,\bm{\Sigma}_N\}$, Agueh and Carlier \cite{agueh2011barycenters} prove that the Wasserstein barycenter exists, is unique, and is also Gaussian with covariance $\bm{\Sigma}_{\text{AC}}$ that uniquely solves
\begin{align}
    \sum_{i=1}^N(\bm{\Sigma}_{\text{AC}}^{\frac{1}{2}}\mathbf{\Sigma}_i^{\vphantom{\frac{1}{2}}}\bm{\Sigma}_{\text{AC}}^{\frac{1}{2}})^{\frac{1}{2}} = \bm{\Sigma}_{\text{AC}}.
\end{align}
For the sake of illustration, we fix $n=2$ and $N=3$. We consider the set of covariances
 \begin{align*}
     \mathbf{\Sigma}_1&= \begin{bmatrix} 0.300 & 0.000 \\ 0.000 & 1.000\end{bmatrix}, & \mathbf{\Sigma}_3&= \begin{bmatrix} \hphantom{-}2.300 & -1.200 \\ -1.200 & \hphantom{-}0.800\end{bmatrix}, \\
     \mathbf{\Sigma}_2&= \begin{bmatrix} 1.600 & 0.600 \\ 0.600 & 0.900\end{bmatrix}.
 \end{align*}
The Agueh-Carlier barycenter is
% and the corresponding optimal transportation maps are
\begin{align*}
    \mathbf{\Sigma}_{\text{AC}}&\approx \begin{bmatrix}\hphantom{-}1.154 & -0.193 \\ -0.193 & \hphantom{-}0.741\end{bmatrix},  \,\,
    % \factor_{\text{ref}\mapsto 1}^\star\approx \begin{bmatrix} 0.519 & 0.081 \\ 0.081 & 1.179\end{bmatrix}, \\
    % \factor_{\text{ref}\mapsto 2}^\star&\approx \begin{bmatrix} 1.202 & 0.403 \\ 0.403 & 1.092\end{bmatrix},  \,\,
    % \factor_{\text{ref}\mapsto 3}^\star\approx \begin{bmatrix} \hphantom{-}1.279 & -0.484 \\ -0.484 & \hphantom{-}0.730\end{bmatrix}.
\end{align*}
With $\Sigma_{\text{ref}} = \bm{\Sigma}_{\text{AC}}$, we proceed to solve Problem \ref{prob:isoholonomic_triangle} for each pair $(i,j)\in\{(1,2),(2,3),(3,1)\}$. We then construct the curve $\Sigma_.$ as in \eqref{eq:wasserstein_polygon_boundary}. 
% To obtain the segment $\Sigma_.^{1\mapsto 2}$ in the curve $\Sigma_.$, we simulate \eqref{eq:IMT_necessary_conditions} with the initial conditions $\Sigma_0=\mathbf{\Sigma}_1$, and 
% \begin{align*}
%     \dot{\Sigma}_0&= \begin{bmatrix} 0.710 & -0.826 \\ -0.826 & -0.525 \end{bmatrix}, & \Omega &= \begin{bmatrix} 0 & -2.214 \\ 2.214 & 0 \end{bmatrix},
% \end{align*}
% over the interval $\tau\in[0,1]$. Similarly, to obtain $\Sigma_.^{2\mapsto 3}$, we use $\Sigma_0=\mathbf{\Sigma}_2$, and 
% \begin{align*}
%     \dot{\Sigma}_0&= \begin{bmatrix} -0.879 & -1.886 \\ -1.886 & -0.141 \end{bmatrix}, & \Omega &= \begin{bmatrix} 0 & -0.690 \\ -0.690 & 0 \end{bmatrix}.
% \end{align*}
% Finally, to obtain $\Sigma_.^{3\mapsto 1}$, we take $\Sigma_0=\mathbf{\Sigma}_3$, and 
% \begin{align*}
%     \dot{\Sigma}_0&= \begin{bmatrix} -2.445 & \hphantom{-}2.072 \\ \hphantom{-}2.072 & -1.005 \end{bmatrix}, & \Omega &= \begin{bmatrix} 0 & -0.917 \\ 0.917 & 0 \end{bmatrix}.
% \end{align*}
The length of the curve $\Sigma_{.}$, i.e. the perimeter of the Wasserstein triangle, is computed to be $\approx 3.567$.

\subsubsection*{The Ando-Li-Mathias geometric mean}

An alternative to the Agueh-Carlier ``averaging'' of distributions, for the special case of Gaussian distributions, is the Ando-Li-Mathias geometric mean \cite{ando2004geometric}. Specifically, for any collection of positive definite matrices $\mathbf{\Sigma}_i$, $i\in\{1,\ldots,N\}$, one can obtain a positive definite matrix $\mathbf{\Sigma}_{ALM}$ with very appealing properties that can be called the $N$-th root of the product of the $\Sigma$'s; when they all commute, it is precisely the $N$-th root of their product and, in general, transforms conformably under congruence transformation.

The construction of $\mathbf{\Sigma}_{ALM}$ is recursive,
starting from two $\mathbf{\Sigma}_i,\mathbf{\Sigma}_j$, in which case it is (see also footnote  \ref{footnote3})
\[
{\rm gm}_2(\mathbf{\Sigma}_i,\mathbf{\Sigma}_j):=
\mathbf{\Sigma}_{i}^{\frac{1}{2}}
    \big(\mathbf{\Sigma}_{i}^{\frac{1}{2}}\mathbf{\Sigma}_j^{-1}\mathbf{\Sigma}_{i}^{\frac{1}{2}}\big)^{-\frac12}
    \mathbf{\Sigma}_{i}^{\frac{1}{2}},
\]
progressing to the geometric mean  ${\rm gm_3}(\cdot,\cdot,\cdot)$ of three matrices, in which case it is the fixed point of the iteration
\[
(\mathbf{\Sigma}_i,\mathbf{\Sigma}_j,\mathbf{\Sigma}_k) \mapsto 
({\rm gm}_2(\mathbf{\Sigma}_j,\mathbf{\Sigma}_k),
{\rm gm}_2(\mathbf{\Sigma}_k,\mathbf{\Sigma}_i),
{\rm gm}_2(\mathbf{\Sigma}_i,\mathbf{\Sigma}_j)),
\]
and similarly to any number $N$ of covariances. For the same set of covariances as in the Agueh-Carlier case, the Andos-Li-Mathias geometric mean is
% and the corresponding optimal transportation maps are
\begin{align*}
    \mathbf{\Sigma}_{\text{ALM}}&\approx \begin{bmatrix}\hphantom{-}0.864 & -0.178 \\ -0.178 & \hphantom{-}0.622\end{bmatrix},  \,\, 
    % \factor_{\text{ref}\mapsto 1}^\star\approx \begin{bmatrix} 0.605 & 0.106 \\ 0.106 & 1.292\end{bmatrix}, \\ 
    % \factor_{\text{ref}\mapsto 2}^\star&\approx \begin{bmatrix} 1.401 & 0.480 \\ 0.480 & 1.207\end{bmatrix}, \,\,
    % \factor_{\text{ref}\mapsto 3}^\star\approx \begin{bmatrix} \hphantom{-}1.462 & -0.534 \\ -0.534 & \hphantom{-}0.803\end{bmatrix}.
\end{align*}
With $\Sigma_{\text{ref}} = \mathbf{\Sigma}_{\text{ALM}}$, we proceed to solve Problem \ref{prob:isoholonomic_triangle} for each pair $(i,j)\in\{(1,2),(2,3),(3,1)\}$. We then construct the curve $\Sigma_.$ as in \eqref{eq:wasserstein_polygon_boundary}. 
% To obtain the segment $\Sigma_.^{1\mapsto 2}$ in the curve $\Sigma_.$, we simulate \eqref{eq:IMT_necessary_conditions} with the initial conditions $\Sigma_0=\mathbf{\Sigma}_1$, and 
% \begin{align*}
%     \dot{\Sigma}_0&= \begin{bmatrix} 0.707 & -0.831 \\ -0.831 & -0.526 \end{bmatrix}, & \Omega &= \begin{bmatrix} 0 & -2.221 \\ 2.221 & 0 \end{bmatrix}.
% \end{align*}
% Similarly, to obtain $\Sigma_.^{2\mapsto 3}$, we use $\Sigma_0=\mathbf{\Sigma}_2$, and 
% \begin{align*}
%     \dot{\Sigma}_0&= \begin{bmatrix} -1.136 & -1.865 \\ -1.865 & -0.049 \end{bmatrix}, & \Omega &= \begin{bmatrix} 0 & -0.981 \\ 0.981 & 0 \end{bmatrix}.
% \end{align*}
% Finally, to obtain $\Sigma_.^{3\mapsto 1}$, we take $\Sigma_0=\mathbf{\Sigma}_3$, and 
% \begin{align*}
%     \dot{\Sigma}_0&= \begin{bmatrix} -2.544 & \hphantom{-}1.974 \\ \hphantom{-}1.974 & -0.853 \end{bmatrix}, & \Omega &= \begin{bmatrix} 0 & -0.715 \\ 0.715 & 0 \end{bmatrix}.
% \end{align*}
The length of the curve $\Sigma_{.}$ in this case is computed as $\approx 3.573$. The isoholonomic Wasserstein triangles corresponding to the two choices of the reference covariance $\Sigma_{\text{ref}}$ are shown in Figure \ref{fig:isoholonomic_triangle_combined} along with a McCann triangle. \footnote{To visualize the trajectory, we use the notation
\begin{align}
    \Sigma=\begin{bmatrix}
        \sigma_1 & \sigma_2 \\ \sigma_2 & \sigma_3
    \end{bmatrix},
\end{align}
to denote the entries of a symmetric matrix.}
% As can be observed, the two isoholonomic triangles are close to one another but are emphatically different from the McCann triangle. The natural next step is to optimize the perimeter of the triangle with respect to $\Sigma_{\text{ref}}$. We pursue this in future work.

\begin{remark}
{ The path along a cycle may be optimized over the parallel transport maps along each segment, subject to the condition of a trivial holonomy. The more restrictive scheme suggested in the current section to achieve a trivial holonomy, i.e., registering the distributions with respect to a reference, simplifies the computation of parallel transport maps.}
     \hfill $\Box$
\end{remark}

\section{Epilogue}\label{sec:Epilogue}
The impact of Monge's problem, to establish an optimal correspondence between distributions, has been foundational across mathematics, engineering, and the sciences \cite{villani2009optimal,peyre2019computational}. It constitutes a basic paradigm in stochastic control \cite{chen2016optimal,chen2021optimal} and has provided a setting for addressing a host of problems in statistics, optimization, probability theory, thermodynamics, economics, machine learning, to name a few.

\tcb{Out of a rapidly growing literature we select two references that highlight the relevance of potential forces, i.e., generated as gradients of a potential, and hence to symmetric $A_t$'s. In the first one, Movilla Miangolarra etal.\ \cite{movilla2021energy} study thermodynamic ensembles driven by potential forces to traverse  a closed thermodynamic cycle, to obtain fundamental bounds on the maximal power that can be drawn from thermal anisotropy. In the next, Lavenant etal.\ \cite{lavenant2024toward} develop a framework of trajectory inference for cell differentiation and development, where the underlying potential forcing is due to the so-called Waddington's landscape. In both situations, the constraint on vector fields governing the evolution of thermodynamic and cell states being potential arises from physical considerations.}
%, see e.g., \cite{terpin2024dynamic,shankar2022optimal,elamvazhuthi2023dynamical,movilla2021energy}.

The motivation of this work has been to explore and bring to life a hitherto unstudied angle of a Lagrangian viewpoint \tcb{in optimal transport}, where the actual trajectories and mixing of \tcb{individual (tracer)} particles are of interest. \tcb{We develop our framework in the simplest setting of linear dynamics and Gaussian distributions.} Thus, our central theme has been to control Gaussian distributions using potential forces, and to characterize the holonomy of particle trajectories. The OMT principal bundle structure duly encodes particle trajectories by recording state transition matrices $\factor_t$. 
Equivalently, the trajectories of $n$ representative tracer particles is encoded in the columns of $\Phi_t$, and these provide information on the holonomy of the flow.

Our main contribution has been to introduce and study the holonomy of the Monge-Kantorovich optimal mass transport. The key insight has been to encode the Eulerian constraint on the velocity field to be potential into the OMT principal Ehresmann connection. We then explored the sub-Riemannian structure induced by this connection along with the Wasserstein-Otto metric in the context of Gaussian distributions. A more comprehensive study of this structure in finite and infinite dimensions, \tcb{as in e.g., \cite{khesin2008geometry},} is of interest, and \tcb{is currently being pursued by the authors. We also note an apparent link between the holonomy of optimal transport and quantum holonomy, as defined by Berry and Simon on the space of density matrices, where a horizontal space of Hermitian vector fields induces a natural Ehresmann connection; see, e.g., \cite[Eq.\ (3)]{uhlmann1986parallel}, which has given rise to growing literature on alternative geometries of positive and invertible elements of unital $C^*$-algebras \cite{corach1999differential}.}

Finally, a similar framework to the one developed herein should prove useful in ensemble control \cite{li2010ensemble,chen2018state,caluya2020finite}, where trajectories of multi-species populations indexed by $s$, for instance, are encoded in collections of transition matrices, $\factor^s_t$'s, on a suitable bundle $\pi$ projecting onto the configuration in a base space, e.g., of covariance data of the form $\Sigma_t=
\sum_s \factor^s_t(\factor^s_t)^\intercal$.
Other variants of practical interest include modeling the flow of particles based on
aggregate information and the location of tracer particles used to estimate the holonomy of the collection, and, likewise, estimation of trajectories of dynamical systems driven by stochastic excitation based on sample trajectories. These are the subject of ongoing work.
\label{seconclusions}

\balance
\bibliographystyle{plain}

\bibliography{References}

\end{document}